\documentclass[preprint,11pt]{elsarticle}

\usepackage{amssymb}
\usepackage{mathrsfs}
\usepackage{mathrsfs}
\usepackage{mathrsfs}
\usepackage{amsfonts}
\usepackage{amsfonts}
\usepackage{amsfonts}
\usepackage{mathrsfs}
\usepackage{amsmath}
\usepackage{graphicx}
\usepackage{multirow}
\usepackage{caption2}
\usepackage{float}
\usepackage{booktabs}
\usepackage{algorithm}
\usepackage{algorithmic}
\usepackage{subfigure}
\newtheorem{theorem}{Theorem}[section]

\numberwithin{equation}{section}

\newtheorem{definition}[theorem]{Definition}

\newtheorem{lemma}[theorem]{Lemma}

\newtheorem{property}[theorem]{Property}

\newenvironment{proof}[1][Proof]{\textbf{#1. }}{\ \rule{0.5em}{0.5em}}%

\journal{}

\begin{document}
\begin{frontmatter}



\title{A Gauss-Seidel Iterative Thresholding Algorithm for $l_q$  Regularized Least Squares Regression
 \tnoteref{t1}} \tnotetext[t1]{The research was supported by the National Natural Science
Foundation of China (Grant No. 11401462).}

\author{Jinshan Zeng$^1$}

\author{Zhiming Peng$^2$}

\author{Shaobo Lin$^3$ \corref{*} }\cortext[*]{Corresponding author:
sblin1983@gmail.com}

\address{1. School of Computer and Information Engineering, Jiangxi Normal University, Nanchang, 330022, P R China.

2. Department of Mathematics, University of California, Los Angeles
(UCLA), Los Angeles, CA 90095, United States.

3.College of Mathematics and Information Science, Wenzhou
University, Wezhou, 325035, P R China }

\begin{abstract}
In  recent studies on sparse modeling, $l_q$ ($0<q<1$) regularized
least squares regression ($l_q$LS) has received considerable
attention due to its superiorities on sparsity-inducing and
bias-reduction over the convex counterparts. In this paper, we
propose a  \underline{Ga}uss-Seidel \underline{i}terative
\underline{t}hresholding \underline{a}lgorithm (called GAITA) for
solution to this problem. Different from the classical iterative
thresholding algorithms using the Jacobi updating rule, GAITA takes
advantage of the Gauss-Seidel rule to update the coordinate
coefficients. Under a mild condition, we can justify that the
support set and sign of an arbitrary sequence generated by GAITA
will converge within finite
  iterations. This convergence property together with the
Kurdyka-{\L}ojasiewicz property of ($l_q$LS) naturally yields the
strong convergence of GAITA  under the same condition as above,
which is generally weaker than the  condition for the convergence of
the classical iterative thresholding algorithms. Furthermore, we
demonstrate that GAITA converges to a local minimizer under certain
additional conditions. A set of numerical experiments are provided
to show the effectiveness, particularly, much faster convergence of
GAITA as compared with the classical iterative thresholding
algorithms.
\end{abstract}

\end{frontmatter}
\section{Introduction}

In this paper, we consider the following $l_q$ ($0<q<1$) regularized
least squares regression ($l_q$LS) problem
\begin{equation}
(l_q\text{LS})\ \ \ \min_{x\in \mathbf{R}^N} \left\{T_{\lambda}(x) =
\frac{1}{2}{\Vert Ax-y\Vert}_{2}^{2} + \lambda \|x\|_q^q \right\},
\label{lqReg}
\end{equation}
where $\|x\|_q^q = \sum_{i=1}^N |x_i|^q$, $N$ is the dimension of
$x$ and $\lambda>0$ is a regularization parameter. The ($l_q$LS)
problem has attracted lots of attention in both scientific research
and engineering practice, since it commonly has stronger
sparsity-promoting ability and better bias-reduction property over
the $l_1$ case. Its typical applications include  signal processing
{\cite{Chartrand2007}, \cite{Chartrand2008}}, image processing
{\cite{L2/3Cao2013}, \cite{Krishnan2009}}, synthetic aperture radar
imaging {\cite{SARZeng2012},
and machine learning
{\cite{Linlq2014}}.

One of the most important class of algorithms to solve the ($l_q$LS)
problem is the iterative thresholding algorithm (ITA)
\cite{BrediesNonconvex2015}, \cite{L1/2TNN}.  Compared with some
other classes of algorithms such as the reweighted least squares
(IRLS) minimization   \cite{Daubechies2010}  and iterative
reweighted $l_1$-minimization (IRL1) {\cite{Candes2008}} algorithms,
ITA generally has lower computational complexity for large scale
problems {\cite{SARZeng2012}}, which triggered   avid research
activities of ITA in the past decade (see {\cite{Blumensath2008,
Daubechies2004Soft, L1/2TNN, ZengHalf2014}}). The  makeup of  ITA
comprises two  steps: a gradient descent-type iteration for  the least
squares and a thresholding operator. To be detailed, for an arbitrary
$\mu>0$, the thresholding function (or proximity operator)
for ($l_q$LS) can be defined as
\begin{equation}
        Prox_{\mu,\lambda \|\cdot\|_q^q}(x) = \arg \min_{u\in \mathbf{R}^N}
        \left\{\frac{\|x-u\|_2^2}{2\mu} + \lambda \|u\|_q^q\right\}.
         \label{ProxOper}
\end{equation}
Since $\|\cdot\|_q^q$ is separable,  computing $Prox_{\mu,\lambda
\|\cdot\|_q^q}$ can be reduced to solve several one-dimensional
minimization problems, that is,
\begin{equation}
          prox_{\mu,\lambda |\cdot|^q}(z) = \arg \min_{v\in \mathbf{R}}
          \left\{\frac{|z-v|^2}{2\mu} + \lambda |v|^q\right\},
\label{SVProxOper}
\end{equation}
and thus,
\begin{equation}
Prox_{\mu,\lambda \|\cdot\|_q^q}(x) = (prox_{\mu,\lambda |\cdot|^q}(x_1), \cdots,
prox_{\mu,\lambda |\cdot|^q}(x_N))^T. \label{ProxOper1}
\end{equation}
For some $q$, such as  $\frac{1}{2}$ or $\frac{2}{3}$,
  $prox_{\mu,\lambda |\cdot|^q}(\cdot)$ can be analytically
expressed  {\cite{L1/2TNN}}. While for other $q\in (0,1)$, we can
use an iterative scheme proposed by {\cite{lqCD-Marjanovic-2014}} to
compute the operator $prox_{\mu,\lambda |\cdot|^q}(\cdot)$. All
these make  the thresholding operator achievable. Then, an efficient
gradient-descent iteration for the un-regularized least squares
problem ($\lambda=0$ in ($l_q$LS)) together with the aforementioned
thresholding operator can derive a feasible scheme to solve
($l_q$LS).
%
%

\subsection{Jacobi iteration and Gauss-Seidel iteration}
As the thresholding operator depends only on $q$, the convergence of
ITA depends heavily on the attributions of the  gradient-descent
type iteration. Landweber-type iteration,   is a natural
selection to solve the un-regularized least squares problems,
since its feasibility has been sufficiently verified in many literatures (say, \cite{Kivinen1997}).
In the classical ITA {\cite{Blumensath2008, Daubechies2004Soft, L1/2TNN}},
a Jacobi iteration strategy whose  Landweber iteration rule is
imposed on the variable $x^n$, is employed to derive the estimate.
 We denote such algorithm as  \textbf{JAITA} henceforth.  More specially, JAITA for
($l_q$LS) can be described as:
\begin{equation}
x^{n+1} \in Prox_{\mu,\lambda \|\cdot\|_q^q}(x^n-\mu A^T(Ax^n-y)), \label{ITAq}
\end{equation}
where $\mu>0$ is a step size parameter.

As a cousin of the Jacobi scheme, the Gauss-Seidel scheme is also
widely used to build blocks for more complicated algorithms
{\cite{Tseng2001, Tseng-Yun2009, Tsitsiklis1989, Xu-Yin2013}.
Different from the Jacobi iteration that updates all the components
simultaneously, the Gauss-Seidel iteration is a component-wise
scheme. Generallly speaking, the Gauss-Seidel iteration is faster than
the corresponding Jacobi iteration {\cite{Tsitsiklis1989}}, since it
uses the latest updates at each iteration.  The aim of this paper is
to introduce the Gauss-Seidel scheme to solve ($l_q$LS). The core
construction of the detailed Gauss-Seidel update rule is by a
concrete representation of the thresholding function, which is
derived by the most recent work {\cite{BrediesNonconvex2015}}.

According to {\cite{BrediesNonconvex2015}}, $prox_{\lambda\mu, q}(\cdot)$ can be expressed as \\
\begin{equation}
prox_{\mu,\lambda |\cdot|^q} (z)=
\left\{
\begin{array}
[c]{ll}%
(\cdot + \lambda \mu q sgn(\cdot) |\cdot|^{q-1})^{-1}(z), & \mbox{for} \ |z|\geq \tau_{\lambda\mu,q}\\
0, & \mbox{for} \ |z|\leq \tau_{\lambda\mu,q}
\end{array}
\right.
\label{ProxMapExpLq}
\end{equation}
for any $z \in \mathbf{R}$
with
\begin{equation}
\tau_{\lambda\mu,q} = \frac{2-q}{2-2q}(2\lambda\mu(1-q))^{\frac{1}{2-q}},
\label{ThreshValuexLq}
\end{equation}
\begin{equation}
\eta_{\lambda\mu,q} = (2\lambda\mu (1-q))^{\frac{1}{2-q}},
\label{ThreshValueyLq}
\end{equation}
and the range of $prox_{\mu,\lambda |\cdot|^q}$ is $\{0\}\cup [\eta_{\lambda\mu,q},\infty)$, $sgn(\cdot)$ represents the sign function henceforth.
When $\ |z|\geq \tau_{\lambda\mu,q}$, the relation $prox_{\mu,\lambda |\cdot|^q} (z)= (\cdot + \lambda \mu q sgn(\cdot) |\cdot|^{q-1})^{-1}(z)$ means that
$prox_{\mu,\lambda |\cdot|^q} (z)$ satisfies the following equation
\[
v + \lambda \mu q\cdot sgn(v) |v|^{q-1} = z.
\]

Now we are in a position to present the proposed algorithm by
utilizing  the Gauss-Seidel iteration. Given the current estimate
$x^n$ and the step size $\mu$, at the next iteration, the $i$-th
coefficient is selected cyclically by
\begin{equation}
  i= \left\{
  \begin{array}{cc}
  N   & {\rm if}\  0\equiv {(n+1)}\  {\rm mod} \ N \\
  {(n+1)}\  {\rm mod} \ N & {\rm otherwise} %
  \end{array}%
  \right..
  \label{updatingindex}
\end{equation}%
We then derive a component-based update of the un-regularized least
squares by
\begin{equation}
z_i^n = x_i^n - \mu A_i^T(Ax^n-y), \label{updateFowStep}
\end{equation}
which together with the thresholding operator then yields a
component-based update for ($l_q$LS) as
\[
x_i^{n+1} \in \arg\min_{v\in \mathbf{R}}
\left\{\frac{|z^n_i-v|^2}{2} + \lambda\mu |v|^q\right\} =
prox_{\mu,\lambda |\cdot|^q}(z^n_i).
\]
  It can be seen from (\ref{ProxMapExpLq}) that
$prox_{\mu,\lambda |\cdot|^q}$ is a set-valued operator. Therefore, motivated
by \cite{lqCD-Marjanovic-2014}, we select a particular single-valued
operator of $prox_{\mu,\lambda |\cdot|^q}$ and then update $x_i^{n+1}$
according to the following scheme,
\begin{equation}
x^{n+1}_i= \mathcal{T} (z_i^n, x_i^n), \label{updatingrul1}
\end{equation}
where
\begin{equation*}
\mathcal{T} (z_i^n, x_i^n)=
  \left\{
  \begin{array}{cc}
  prox_{\mu,\lambda |\cdot|^q}(z_i^{n})   & {\rm if} \ |z_i^{n}| \neq \tau_{\lambda\mu,q} \\
  sgn(z_i^{n})\eta_{\lambda\mu,q} \mathbf{I}(x_i^{n} \neq 0), & {\rm if} \ |z_i^{n}| = \tau_{\lambda\mu,q}%
  \end{array}%
  \right.,
\end{equation*}%
and $\mathbf{I}(x_i^{n} \neq 0)$ denotes the indicator function,
that is,
\begin{equation*}
  \mathbf{I}(x_i^{n} \neq 0)= \left\{
  \begin{array}{cc}
  1,   & {\rm if}\  x_i^{n} \neq 0\\
  0, & {\rm otherwise} %
  \end{array}%
  \right..
\end{equation*}%
While the other components of $x^{n+1}$ are fixed, i.e.,
\begin{equation}
x_j^{n+1} = x_j^{n}, \ {\rm for}\ j\neq i. \label{updatingrul2}
\end{equation}
For the sake of brevity, we denote in the rest of paper
$\tau_{\mu,q}$ and $\eta_{\mu,q}$ to take the place of
$\tau_{\lambda\mu,q}$ and $\eta_{\lambda\mu,q}$, respectively. In
summary, we can formulate the proposed algorithm as follows.

\begin{center}
\text{\underline{Ga}uss-Seidel \underline{I}terative
\underline{T}hresholding \underline{A}lgorithm (GAITA)}\\
\vspace{0.1cm}
\begin{tabular}{l}
  \hline
  Initialize with $x^0$. Choose a step size $\mu>0$, let $n:=0$. \\ \vspace{0.2cm}
  Step 1. Calculate the index $i$ according to (\ref{updatingindex}); \\ \vspace{0.2cm}
  Step 2. Calculate $z_i^n$ according to (\ref{updateFowStep}); \\ \vspace{0.2cm}
  Step 3. Update $x_i^{n+1}$ via (\ref{updatingrul1}) and $x_j^{n+1} = x_j^{n}$  for $j\neq i$;\\ \vspace{0.2cm}
  Step 4. Check the terminational rule. If yes, stop; \\ \vspace{0.2cm}
  \ \ \ \ \ \ \ \ \ otherwise, let $n:= n+1$, go to Step 1.\\
  \hline
\end{tabular}\\
\end{center}

\subsection{Why Gauss-Seidel?}
It can be found in the last section that the main difference between
JAITA and GAITA lies at  whether the Landweber iteration is
component-wise. Such a slight difference  leads to a plausible
assertion that the convergence of both algorithms are similar. To
verify the authenticity of the above viewpoint, we conduct a set of
experiments to    the convergence of JAITA and GAITA. Interestingly,
we find in this experiment that the convergence of the
aforementioned two algorithms are totally different.

 To be
detailed, given a sparse signal $x$ with dimension $N=500$ and
sparsity $k^*=15,$ we considered  the signal recovery problem through
the observation $y=Ax,$ where the original sparse signal $x$ was
generated randomly according to the standard Gaussian distribution,
and $A$ was of dimension $m\times N=250 \times 500$ with Gaussian
$\mathcal{N}(0,1/250)$ i.i.d. entries and was preprocessed via
column-normalization, i.e., $\|A_i\|_2 =1$ for any $i$. We then
applied GAITA and JAITA  to the ($l_q$LS) problem with two different
$q$, that is, $q=1/2$ and $2/3$, respectively. In both cases, the
thresholding functions can be analytically expressed as shown in
{\cite{L1/2TNN}} and {\cite{L2/3Cao2013}}, respectively, and thus
the corresponding algorithms can be efficiently implemented. In both
cases, we set $\lambda=0.001$, $\mu = 0.95$ for both JAITA and
GAITA. Moreover, the initial guesses were taken as 0 for all cases.
The trends of the objective sequences in different cases are shown
in Fig. {\ref{Weaker_Cond}}.

From Fig. {\ref{Weaker_Cond}}, the objective sequences of JAITA
diverge for both $q=1/2$ and $2/3$, while those of GAITA are
definitely convergent. This means that there exists some  $\mu$ such
that JAITA is divergent but GAITA is assuredly convergent, which
significantly  stimulates our research interests, since a large
scope of $\mu$ to guarantee the convergence essentially enlarges
the applicable range of iterative thresholding-type algorithms.




We then naturally turn to theoretically  verify the interesting
phenomenon shown by Fig.{\ref{Weaker_Cond}}.  That is, the aim of
our study is to answer the following questions:

\begin{enumerate}
\item[(Q1)]
Is the convergence condition of GAITA exactly  weaker than that of
JAITA?
\end{enumerate}

\begin{enumerate}
\item[(Q2)]
If the answer of the above question is positive, then what is the
applicable range of $\mu$ for GAITA to guarantee the convergence?
\end{enumerate}


\begin{figure}[!t]
\begin{minipage}[b]{.49\linewidth}
\centering
\includegraphics*[scale=0.45]{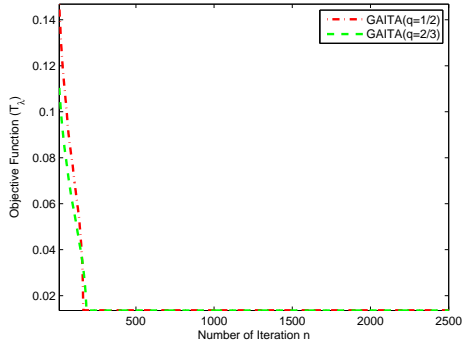}
\centerline{{\small (a) Convergence of GAITA}}
\end{minipage}
\hfill
\begin{minipage}[b]{.49\linewidth}
\centering
\includegraphics*[scale=0.45]{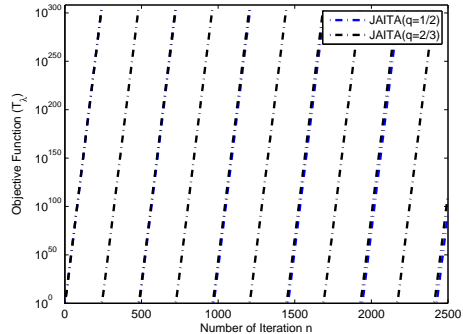}
\centerline{{\small (b) Divergence of JAITA}}
\end{minipage}
\hfill \caption{ An experiment that motivates the use of
the Gauss-Seidel scheme.
(a) The trends of the objective function sequences, i.e.,
$\{T_{\lambda}(x^{n})\}$ of GAITA for different $q$. (b) The trends
of the objective function sequences of JAITA for different $q$.
} \label{Weaker_Cond}
\end{figure}

\subsection{Related Literatures}

There are many methods  used to solve the ($l_q$LS) problem. Some
general methods such as those in {\cite{Attouch2013, Bagirov2013,
Burke2005, BrediesNonconvex2015, Candes2008, Chen2012,
Daubechies2010, Fuduli2004}} and references therein and also books
{\cite{Bertsekas1999, Nocedal-Wright2006}} do not update the
iterations by using the Gauss-Seidel scheme. In
{\cite{BrediesNonconvex2015}}, the subsequential convergence of the
iterative thresholding algorithm for ($l_q$LS) with an arbitrary
$q\in (0,1)$ and further the global convergence for ($l_q$LS) with a
rational $q$ have been verified under the condition
$0<\mu<\|A\|_2^{-2}$. In {\cite{Attouch2013}}, the global
convergence of the iterative thresholding algorithm for ($l_q$LS)
with an arbitrary $q$ has been justified under the same condition.
Besides these general methods, there are several specific iterative
thresholding algorithms for solving ($l_q$LS) with a specific $q$
such as {\em hard} for $l_0$ {\cite{Blumensath2008}}, {\em soft} for
$l_1$ {\cite{Daubechies2004Soft}} and {\em half} for $l_{1/2}$
{\cite{L1/2TNN}}. Under the same condition, all these specific
iterative thresholding algorithms converge to a stationary point.

Another tightly related class of algorithms is the block coordinate
descent (BCD) algorithm. BCD has been numerously used in many
applications. Its original form, block coordinate minimization (BCM)
can date back to the 1950's {\cite{Hildreth1957}}. The main idea of
BCM is to update a block by minimizing the original objective with
respect to that block. Its convergence was extensively studied under
many different cases (cf. {\cite{Grippo1999},
\cite{Razaviyayn2013}}, \cite{Tseng2001}, {\cite{Xu-Yin2013}} and
the references therein). In {\cite{Luo-Tseng1992}}, the convergence
rate of BCM was developed under the strong convexity assumption for
the multi-block case, and in {\cite{Beck2013}}, its convergence rate
was established under the general convexity assumption for the
two-block case.
Besides BCM, the block coordinate gradient descent (BCGD) method was
also largely studied (cf. {\cite{Tseng-Yun2009}}).
Different from BCM, BCGD updates a block via taking a block gradient
step, which is equivalent to minimizing a certain prox-linear
approximation of the objective. Its global convergence was justified
under the assumptions of the so-called local Lipschitzian error
bound and the convexity of the non-differentiable part of the
objective. In {\cite{Nesterov2012}}, a randomized block coordinate
descent (RBCD) method was proposed. RBCD randomly chooses the block
to update with positive probability at each iteration and is not
essentially cyclic. The weak convergence was established in
{\cite{Nesterov2012}}, {\cite{Richtarik2012}}, while there is no
strong convergence result for RBCD.

One important subclass of BCD is the cyclic coordinate descent (CCD)
algorithm. The CCD algorithm updates the iterations by the cyclic
coordinate updating rule. The work {\cite{Xu-Yin2013}} used cyclic
updates of a fixed order and supposes block-wise convexity. In
{\cite{ncxCDMazumder2007}}, a CCD algorithm was proposed for a class
of non-convex penalized least squares problems. However, both
{\cite{ncxCDMazumder2007}} and {\cite{Tseng2001}} did not consider
the CCD algorithm for the ($l_q$LS) problem. In
{\cite{l1CDFriedman2007}}, a CCD algorithm was implemented for
solving the ($l_1$LS) problem. Its convergence can be shown by
referring to {\cite{Tseng2001}}. In {\cite{Seneviratne2013}}, the
$l_0$LS-CD algorithm was proposed for the ($l_0$LS) problem, and its
convergence to a local minimizer was also shown under certain
conditions. Recently, Marjanovic and Solo
{\cite{lqCD-Marjanovic-2014}} proposed a cyclic descent algorithm
(called $l_q$CD) for the ($l_q$LS) problem with $0<q<1$ and $A$
being column-normalized, i.e., $\|A_i\|_2 = 1$, $i=1,2,\dots, N$,
where $A_i$ is the $i$-th column of $A$. They proved the
subsequential convergence and further the convergence to a local
minimizer under the so-called scalable restricted isometry property (SRIP)
in {\cite{lqCD-Marjanovic-2014}}. In the perspective of the
iterative form, $l_q$CD is a special case of GAITA with $A$ being
column-normalized and $\mu=1.$


\subsection{Contributions}

The main contribution of this paper is to present the convergence
analysis of   GAITA for solving the ($l_q$LS) problem.
The finite step convergence of the support set and sign can be
verified under the condition that the step size $\mu$ is less than
$\frac{1}{\max_i \|A_i\|_2^2}$ (see Theorem
{\ref{Thm_SupportConv}}). It means that the support sets and signs
of the sequence $\{x^n\}$ generated by GAITA certainly converge and
remain the same within the finite iterations. Such property is very
important since it can bring a possible way to construct an
auxiliary sequence, which lies in a special subspace and has the
same convergence behavior of the original sequence $\{x^n\}$. Then
with the help of the  Kurdyka-{\L}ojasiewicz property (See Appendix
G) of   $T_{\lambda}$, we can verify the global convergence of GAITA
  under the same condition, i.e.,
$0<\mu<\frac{1}{\max_i \|A_i\|_2^2}$ (See Theorem
{\ref{Thm_MainConv}}).
It can be noted that this condition is generally weaker than that
  of JAITA  (i.e., $0<\mu<\|A\|_2^{-2}$)
{\cite{Attouch2013}}. This gives positive answers to question (Q1)
and (Q2). The improvement on the convergence condition is commonly
very important. It may improve not only the rate-of-convergence but
also the applicability of GAITA as compared with JAITA. Furthermore,
we can also justify that the proposed algorithm converges to a local
minimizer under certain a second-order condition (See Theorem
{\ref{Thm_LocalMin}}). More specifically, let $x^*$ be the limit
point and $I$ be its support set. Then the condition can be
described as: $A_I^TA_I + \lambda q (q-1) \Lambda(x_I^*)$ is
positive definite, where $A_I$ represents the submatrix of $A$ with
column restricted to the index set $I$, $x_I^*$ is the subvector of
$x^*$ restricted to $I$, and $\Lambda(x_I^*)$ is a diagonal matrix
with $(|x_i^*|^{q-2})_{i\in I}$ as the diagonal vector. Besides this
condition, we also give another two sufficient conditions to
guarantee that the limit point is a local minimizer. The
effectiveness, particularly, the faster convergence and weaker
convergence condition of GAITA than JAITA have also been
demonstrated by a series of numerical experiments. All these results
show that utilizing the Gauss-Seidel iteration in ITA for solving
($l_q$LS) is feasible and efficient.

\subsection{Organization}

The remainder of this paper is organized as follows. Some
preliminaries are given in section 2.
In section 3, we give the convergence analysis of GAITA. In section
4, a series of simulations are implemented to demonstrate the
effectiveness of the proposed algorithm. We conclude this paper in
section 5, and present all proofs in Appendix.


\section{Preliminaries}

In this section, we present some preliminaries, which serve as the
basis of the convergence analysis in the next section.

With the definition of the thresholding function (\ref{ProxOper}),
we can define a new operator $G_{\mu,\lambda \|\cdot\|_q^q}(\cdot)$ as
\begin{equation}
G_{\mu,\lambda \|\cdot\|_q^q}(x) = Prox_{\mu,\lambda \|\cdot\|_q^q}(x-\mu A^T(Ax-y))
\label{G_oper}
\end{equation}
for any $x\in \mathbf{R}^N$. We denote ${\mathcal F}_q$ as the fixed
point set of the operator $G_{\mu,\lambda \|\cdot\|_q^q}$, i.e.,
\begin{equation}
{\mathcal F}_q = \{x: x = G_{\mu,\lambda \|\cdot\|_q^q}(x)\}.
\label{FixedPointSet}
\end{equation}


By the definition of $Prox_{\mu,\lambda \|\cdot\|_q^q}$, a type of optimality
conditions of the ($l_q$LS) problem has been derived in
{\cite{lqCD-Marjanovic-2014}}.

\begin{lemma} \label{Lemma_OptimalCondition}
(Theorem 3 in {\cite{lqCD-Marjanovic-2014}}). Given a point $x^*$,
define the support set of $x^*$ as $Supp(x^*) = \{i: x_i^* \neq
0\}$, then $x^* \in {\mathcal F}_q$ if and only if the following
three conditions hold.
\begin{enumerate}
\item[(a)]
For $i\in Supp(x^*)$, $|x_i^*| \geq \eta_{\mu,q}$.

\item[(b)]
For $i\in Supp(x^*)$, $A_i^T(Ax^*-y) + \lambda q sgn(x_i^*)
|x_i^*|^{q-1} = 0 $.

\item[(c)]
For $i\in Supp(x^*)^c$, $|A_i^T(Ax^*-y)| \leq \tau_{\mu,q}/\mu$.
\end{enumerate}
\end{lemma}

We call $x^*$ a {\bf stationary point} of the ($l_q$LS) problem
henceforth if it satisfies the optimality conditions in Lemma
{\ref{Lemma_OptimalCondition}}. Similarly, according to the
definition of the operator $prox_{\mu,\lambda |\cdot|^q}(\cdot)$,
(\ref{ProxMapExpLq}), and the updating rule of GAITA
(\ref{updatingindex})-(\ref{updatingrul2}), we can claim that
$x^{n+1}$ satisfies the following property.

\begin{property}\label{Prop_OptimalCondition_n+1}
Given the current iterate $x^n$ ($n \in \mathbf{N}$), the index set
$i$ is determined via (\ref{updatingindex}), then $x_i^{n+1}$
satisfies either
\begin{enumerate}
\item[(a)]
$x_i^{n+1} = 0,$ or,
\item[(b)]
$|x_i^{n+1}| \geq \eta_{\mu,q}$ and also satisfies the following
equation
\begin{align}
&A_i^T(Ax^{n+1}-y) + \lambda q sgn(x_i^{n+1}) |x_i^{n+1}|^{q-1} \nonumber\\
&= (\frac{1}{\mu}-A_i^TA_i)(x_i^n - x_i^{n+1}). \label{OptCon_i}
\end{align}
that is, $\nabla_i T_{\lambda}(x^{n+1}) =
(\frac{1}{\mu}-A_i^TA_i)(x_i^n - x_i^{n+1}),$ where $\nabla_i
T_{\lambda}(x^{n+1})$ represents the gradient of $T_{\lambda}$ with
respect to the $i$-th coordinate at the point $x^{n+1}$.
\end{enumerate}
\end{property}
As shown by Property \ref{Prop_OptimalCondition_n+1}, the
coordinate-wise gradient of $T_{\lambda}$ with respect to the $i$-th
coordinate at $x^{n+1}$ is not exact zero but with a relative error.
This property can be easily derived from the definition of
$prox_{\mu,\lambda |\cdot|^q}(\cdot)$ and the specific iterative form of
GAITA. More specifically, according to (\ref{ProxMapExpLq}) and
(\ref{updatingrul1}), it holds obviously either $x_i^{n+1} =0$ or
$|x_i^{n+1}| \geq \eta_{\mu,q}$. Moreover, when $|x_i^{n+1}| \geq
\eta_{\mu,q}$, according to (\ref{SVProxOper}), $x_i^{n+1}$ is a
minimizer of the optimization problem (\ref{SVProxOper}) with
$z=z_i^n$. Therefore, $x_i^{n+1}$ should satisfy the following
optimality condition
\begin{equation}
x_i^{n+1} - z_i^n  + \lambda \mu q sgn(x_i^{n+1}) |x_i^{n+1}|^{q-1}
= 0. \label{OptCon1}
\end{equation}
Plugging (\ref{updateFowStep}) into (\ref{OptCon1}) gives
\begin{equation}\label{OptCon2}
       A_i^T(Ax^{n+1}-y) + \lambda q sgn(x_i^{n+1}) |x_i^{n+1}|^{q-1}
       = \frac{1}{\mu}(x_i^n - x_i^{n+1})-A_i^TA(x^n-x^{n+1}).
\end{equation}
Combining (\ref{updatingrul2}) and (\ref{OptCon2}) implies
(\ref{OptCon_i}).

\section{Convergence Analysis}


In this section, we first show the subsequential convergence of
GAITA, then prove its global convergence, and further justify that
the algorithm can converge to a local minimizer.

\subsection{Subsequential Convergence}

To aid the description, we  show that the sequence
$\{T_{\lambda}(x^n)\}$ satisfies the sufficient decrease property
{\cite{Burke1985SuffCon}} at first.

\begin{property}\label{Prop_SufficientDecrease}
Let $\{x^n\}$ be a sequence generated by GAITA. Assume that
$0<\mu<L_{\max}^{-1}$, then
\[
T_{\lambda}(x^n) - T_{\lambda}(x^{n+1}) \geq
\frac{1}{2}(\frac{1}{\mu} - L_{\max}) \|x^n - x^{n+1}\|_2^2,
~\forall n \in \mathbf{N},
\]
where $L_{\max} = \max_{i} \|A_i\|_2^2.$
\end{property}

The proof of this property is presented in Appendix \ref{APPENDIXA}.
From this property, we can claim that the objective sequence
$\{T_{\lambda}(x^n)\}$ converges since it is lower bounded by 0,
that is,   GAITA is weakly convergent. Furthermore, if the
initialization of the sequence is bounded, then based on Property
{\ref{Prop_SufficientDecrease}}, it can easily derive the following
boundedness and asymptotically regular properties of the sequence.

\begin{property}\label{Prop_AsymptoticRegular}
Let $\{x^n\}$ be a sequence generated by GAITA with a bounded
initialization. Assume $0<\mu<L_{\max}^{-1}$, then $\{x^n\}$ is
bounded for any $n\in \mathbf{N}$, and
$$
\sum_{k=0}^n \|x^{k+1} - x^{k}\|_2^2 \leq \frac{2\mu}{1-\mu
L_{\max}} T_{\lambda}(x^0),
$$
and also
$$
\|x^n - x^{n+1}\|_2 \rightarrow 0, \ {\rm as}\ n\rightarrow +\infty.
$$
\end{property}
%
%


The boundedness of $\{x^n\}$ is mainly due to the sufficient
decrease property, the coercivity of $T_{\lambda}$ and the
boundedness assumption of the initialization. While the asymptotic
regular property is mainly due to the sufficient decrease property
and the boundedness of the initialization. From Properties
{\ref{Prop_SufficientDecrease}} and {\ref{Prop_AsymptoticRegular}},
we can justify  the subsequential convergence of GAITA.

\begin{theorem}\label{Thm_SubsequentialConv}
Let $\{x^n\}$ be a sequence generated by GAITA with a bounded
initialization. Assume that $0<\mu<L_{\max}^{-1}$, then the sequence
$\{x^n\}$ has a convergent subsequence. Moreover, let $\mathcal{L}$
be the set of the limit points of $\{x^n\}$, then $\mathcal{L}$ is
closed and connected.
\end{theorem}

The proof of this theorem is presented in Appendix \ref{APPENDIXB}.
This theorem only shows the subsequential convergence of GAITA.
Moreover, we note that $\mathcal{L}$ might not be a set of isolated
points. Due to this,    it  becomes  challenging to justify the
global convergence of GAITA \cite{ZengIJT2014}. More specifically,
there are still two questions on the convergence of the proposed
algorithm:
\begin{enumerate}
\item[(a)]
When does the algorithm converge globally? Under what conditions,
GAITA converges strongly
  in the sense that the whole sequence generated, regardless of the initial point, is convergent.

\item[(b)]
Where does the algorithm converge? Does the algorithm converge to a
global minimizer or more practically, a local minimizer due to the
non-convexity of the optimization problem?
\end{enumerate}

\subsection{Global Convergence}

In this subsection, we will focus on answering the first  question
proposed in the end of the last subsection. More specifically, we
will show that the whole sequence $\{x^n\}$ generated by GAITA
converges as long as the step size $\mu \in (0, L_{\max}^{-1})$.

Given the current iteration $x^n$, we define the descent function as
\begin{equation}
\Delta(x^n, x^{n+1}) = T_{\lambda}(x^n) - T_{\lambda}(x^{n+1}).
\label{DescentFun}
\end{equation}
Note that $x^n$ and $x^{n+1}$ differ only in their $i$-th
coefficient which is determined by (\ref{updatingindex}).
From now on, if not stated, it is assumed that $x_i^{n+1}$ is given
by (\ref{updatingrul1}) and $i$ is given by (\ref{updatingindex}).
The following lemma presents an important property of the descent
function.

\begin{lemma}\label{Lemma_DescentFun}
Let $\{x^n\}$ be a sequence generated by GAITA. Assume that
$0<\mu<L_{\max}^{-1}$, then
\[
\Delta(x^n, x^{n+1}) = 0 \ {\rm if~ and~ only~ if} \ x_i^{n+1} =
x_i^n.
\]
\end{lemma}
The proof of this lemma is obvious. On one hand, if $x_i^{n+1} =
x_i^n$, then $x^{n+1} = x^n$, and thus $\Delta(x^n, x^{n+1}) = 0$.
On the other hand, if $\Delta(x^n, x^{n+1}) = 0$, then Property
{\ref{Prop_SufficientDecrease}} implies $x^{n+1} = x^n$ and thus,
$x_i^{n+1} = x_i^n$.

Moreover, similar to Theorem 10 in {\cite{lqCD-Marjanovic-2014}}, we
can claim that the mapping $\mathcal{T}(\cdot,\cdot)$ is a closed
mapping, shown as follows.

\begin{lemma}\label{Lemma_ClosedMap}
$\mathcal{T}(\cdot,\cdot)$ is a closed mapping, i.e., assume
\begin{enumerate}
\item[(a)]
$x_i^n \rightarrow x_i^*$ as $n \rightarrow \infty;$

\item[(b)]
$x_i^{n+1} \rightarrow x_i^{**}$ as $n\rightarrow \infty$, where
$x_i^{n+1} = \mathcal{T}(z_i^n, x_i^n).$
\end{enumerate}
Then $x_i^{**} = \mathcal{T}(z_i^*, x_i^*)$, where $z_i^* = x_i^* -
\mu A_i^T(Ax^*-y)$.
\end{lemma}

The proof is the essentially the same as that of Theorem 10 in
{\cite{lqCD-Marjanovic-2014}}. The only difference is that
$prox_{\mu,\lambda |\cdot|^q}$ is discontinuous at $\tau_{\mu,q}$ while
$prox_{1,\lambda |\cdot|^q}$ is discontinuous at $\tau_{1,q}$. Therefore, the
closedness of the operator $\mathcal{T}(\cdot,\cdot)$ can not be
changed after introducing a stepsize $\mu$. The following theorem
shows that any limit point of the sequence $\{x^n\}$ is a stationary
point of the ($l_q$LS) problem.

\begin{theorem}\label{Thm_LimitPoint_StatPoint}
Let $\{x^n\}$ be a sequence generated by GAITA with a bounded
initialization, and $\mathcal{L}$ be its limit point set. Suppose
that $0<\mu<L_{\max}^{-1}$, then $\mathcal{L} \subseteq
{\mathcal{F}}_q$.
\end{theorem}

The proof of this theorem is similar to that of Theorem 5 in
{\cite{lqCD-Marjanovic-2014}}. For the completion, we provide the
proof in Appendix \ref{APPENDIXC}.

In the following theorem, we justify the finite step convergence of
the support sets and signs of the sequence $\{x^n\}$, that is, the
support sets and signs of $\{x^n\}$ will converge and remain the
same within a finite iterations.

\begin{theorem}\label{Thm_SupportConv}
Let $\{x^n\}$ be a sequence generated by GAITA with a bounded
initialization. Assume that $0<\mu< L_{\max}^{-1}$ and $x^*$ is any
limit point of $\{x^n\}$, then there exists a sufficiently large
positive integer $n^*>N$ such that when $n > n^*$, it holds
\begin{enumerate}
\item[(a)]
either $x_j^n=0$ or $|x_j^n|\geq \eta_{\mu,q}$ for $j=1,2,\cdots,N;$

\item[(b)]
$I^n=I$;

\item[(c)]
$sgn(x^n) = sgn(x^*)$,
\end{enumerate}
where $I^n = Supp(x^n) = \{i: |x_i^n| \neq 0, i=1,2\cdots, N\}$ and
$I = Supp(x^*)$.
\end{theorem}

The proof of this theorem is shown in Appendix \ref{APPENDIXD}. Form
this theorem, it can be observed that when $n$ is sufficiently
large, the generated sequence $\{x^n\}$ as well as its limit points
will lie in the same subspace $S\subset \mathbf{R}^N$, which has
some special structure. Due to this, it brings a possible way to
construct an auxiliary sequence that has the same convergence
behavior of the original sequence $\{x^n\}$. Thus, we only need to
verify the convergence of the constructed auxiliary sequence instead
of $\{x^n\}$.
The construction of the  auxiliary sequence is a bit standard and is
motivated by \cite{ZengIJT2014}. To be detailed, the sequence can be
constructed according to the following procedure.
\begin{enumerate}
\item[(a)]
Let $n_0 = j_0N>n^*$ for some positive integer $j_0$. Then we can
define a new sequence $\{\hat{x}^n\}$ with $\hat{x}^n = x^{n_0+n}$
for ${n\in \mathbf{N}}$. It is obvious that $\{\hat{x}^n\}$ has the
same convergence behavior with $\{x^n\}$. Moreover, it can be noted
from Theorem {\ref{Thm_SupportConv}} that all the support sets and
signs of $\{\hat{x}^n\}$ are the same.

\item[(b)] Denote $I$ as the convergent support set of the sequence $\{x^n\}$. Let $K$ be the number of elements of $I$.
Without loss of generality, we assume
\[
1\leq I(1)<I(2)<\cdots < I(K) \leq N.
\]
According to the updating rule
(\ref{updatingindex})-(\ref{updatingrul2}) of GAITA, we can observe
that many successive iterations of $\{\hat{x}^n\}$ are the same.
Thus, we can merge these successive iterations into a single
iteration. Moreover, the updating rule of the index is cyclic and
thus periodic. As a consequence, the merging procedure can be
repeated periodically. Formally, we consider such a periodic
subsequence with $N$-length of $\{\hat{x}^n\}$, i.e.,
$$\{\hat{x}^{jN+I(1)},\hat{x}^{jN+I(1)+1}, \cdots, \hat{x}^{jN+I(1)+N-1}\}$$
for $j\in \mathbf{N}$. Then for any $j\in \mathbf{N}$, we emerge the
$N$-length sequence $\{\hat{x}^{jN+I(1)}, \cdots,
\hat{x}^{jN+I(1)+N-1}\}$ into a new $K$-length sequence
$\{\bar{x}^{jK+1},\bar{x}^{jK+2}, \cdots, \bar{x}^{jK+K}\}$ with the
rule
$$
\{\hat{x}^{jN+I(i)},\cdots, \hat{x}^{jN+I(i+1)-1}\} \mapsto
\bar{x}^{jK+i},
$$
with $\bar{x}^{jK+i}=\hat{x}^{jN+I(i)}$ for $i=1,2,\cdots,K,$ since
$\hat{x}^{jN+I(i)+k} = \hat{x}^{jN+I(i)}$ for
$k=1,\cdots,I(i+1)-I(i)-1.$ Moreover, we emerge the first $I(1)$
iterations of $\{\hat{x}^n\}$ into $\bar{x}^0$, i.e.,
$$
\{\hat{x}^{0},\cdots, \hat{x}^{I(1)-1}\} \mapsto \bar{x}^{0},
$$
with $\bar{x}^{0} = \hat{x}^{0}$, since these iterations keep
invariant and are equal to $\hat{x}^{0}$. After this procedure, we
obtain a new sequence $\{\bar{x}^{n}\}$ with $n=jK+i$,
$i=0,\cdots,K-1$ and $j\in \mathbf{N}$. It can be observed that such
an emerging procedure keeps the convergence behavior of
$\{\bar{x}^n\}$ the same as that of $\{\hat{x}^n\}$ and $\{x^n\}$.

\item[(c)] Furthermore, for the index set $I$, we define a projection $P_I$ as
$$
P_I: \mathbf{R^N} \rightarrow \mathbf{R}^K, P_I x = x_I, \forall
x\in \mathbf{R}^N,
$$
where $x_I$ represents the subvector of $x$ restricted to the index
set $I$. With this projection, a new sequence $\{u^n\}$  is
constructed such that
$$
u^n = P_I \bar{x}^n,
$$
for $n\in \mathbf{N}$. As we can observe that $u^n$ keeps all the
non-zero elements of $\bar{x}^n$ while gets rid of its zero
elements. Moreover, this operation can not change the convergence
behavior of $\{\bar{x}^n\}$ and $\{u^n\}$. Therefore, the
convergence behavior of $\{u^n\}$ is the same as $\{x^n\}$.
\end{enumerate}

After the construction procedure (a)-(c), we get a new sequence
$\{u^n\}$. In the following, we will prove the convergence of
$\{x^n\}$ via justifying the convergence of $\{u^n\}$. Let
$$\mathcal{U} = \{u^*: u^* = P_I x^*, \forall x^* \in \mathcal{L}\}.$$
Then $\mathcal{U}$ is the corresponding limit point set of
$\{u^n\}$. Furthermore, we define a new function $T$ as follows:
\begin{equation}
T: \mathbf{R}^K \rightarrow \mathbf{R}, T(u) = T_{\lambda}(P_I^T u),
\forall u\in \mathbf{R}^K, \label{def-T}
\end{equation}
where $P_I^T$ denotes the transpose of the projection $P_I$, and is
defined as
\[
P_I^T: \mathbf{R}^K \rightarrow \mathbf{R}^N, (P_I^T u)_I = u,
(P_I^T u)_{I^c} = 0, \forall u \in \mathbf{R}^K.
\]
Here $I^c$ represents the complementary set of $I$, i.e., $I^c =
\{1,2,\cdots, N\} \setminus I$, $(P_I^T u)_I$  and $(P_I^T u)_{I^c}$
represent the subvectors of $P_I^T u$ restricted to $I$ and $I^c$,
respectively. Let $B=A_I$, where $A_I$ denotes the submatrix of $A$
restricted to the index set $I$. Thus,
\[
T(u) = \frac{1}{2} \|Bu-y\|_2^2 + \lambda \|u\|_q^q.
\]

After the construction procedure (a)-(c), we can observe that the
following properties still hold for $\{u^n\}$.

\begin{lemma}\label{Lemma_propofnewsequence}
The sequence $\{u^n\}$ possesses the following properties:
\begin{enumerate}
\item[(a)]
$\{u^n\}$ is updated via the following cyclic rule. Given the
current iteration $u^n$, only the $i$-th coordinate will be updated
while the other coordinate coefficients will be fixed at the next
iteration, i.e.,
\begin{equation}
u^{n+1}_i= \mathcal{T} (v_i^n, u_i^n), \label{newupdatingrul1}
\end{equation}
and
\begin{equation}
u_j^{n+1} = u_j^{n}, \ {\rm for}\ j\neq i, \label{newupdatingrul2}
\end{equation}
where $i$ is determined by
\begin{equation}
  i= \left\{
  \begin{array}{cc}
  K   & {\rm if}\  0\equiv {(n+1)}\  {\rm mod} \ K \\
  {(n+1)}\  {\rm mod} \ K, & {\rm otherwise} %
  \end{array}%
  \right.,
  \label{newupdatingindex}
\end{equation}%
and
\begin{equation}
v_i^n = u_i^n - \mu B_i^T(Bu^n-y), \label{newupdateFowStep}
\end{equation}

\item[(b)]
According to the updating rules
(\ref{newupdatingrul1})-(\ref{newupdateFowStep}), for $n\geq K$,
  there exit two positive integers $1\leq i_0 \leq K$ and $j_0
\geq 1$ such that $n=j_0K+i_0$ and
\begin{equation}
  u_j^n= \left\{
  \begin{array}{cc}
  u_j^{n-(i_0-j)},   & {\rm if}\  1\leq j \leq i_0\\
  u_j^{n-K-(i_0-j)}, & {\rm if}\ i_0+1 \leq j \leq K %
  \end{array}%
  \right..
  \label{cyclic-u}
\end{equation}%

\item[(c)]
For any $n \in \mathbf{N}$,
\[
u^n \in \mathbf{R}_{{\eta_{\mu,q}}^c}^K,
\]
where $\mathbf{R}_{{\eta_{\mu,q}}^c}$ represents a one-dimensional
real subspace, which is defined as
\[
\mathbf{R}_{{\eta_{\mu,q}}^c} = \mathbf{R} \setminus (-\eta_{\mu,q},
\eta_{\mu,q}).
\]

\item[(d)] Given $u^n$, if $i$ is determined by (\ref{newupdatingindex}), then $u_i^{n+1}$
satisfies the following equation
\begin{align}
&B_i^T(Bu^{n+1}-y) + \lambda q sgn(u_i^{n+1}) |u_i^{n+1}|^{q-1} \nonumber\\
&= (\frac{1}{\mu}-B_i^TB_i)(u_i^n - u_i^{n+1}). \label{newOptCon_i}
\end{align}
That is,
$$
\nabla_i T(u^{n+1}) = (\frac{1}{\mu}-B_i^TB_i)(u_i^n - u_i^{n+1}),
$$
where $\nabla_i T(u^{n+1})$ represents the gradient of $T(\cdot)$
with respect to the $i$-th coordinate at the point $u^{n+1}$.

\item[(e)]
$\{u^n\}$ satisfies the following sufficient decrease condition:
\[
T(u^n) - T(u^{n+1}) \geq a \|u^n - u^{n+1}\|_2^2,
\]
for $n\in \mathbf{N}$, where $a=\frac{1}{2}(\frac{1}{\mu} -
L_{\max}).$

\item[(f)]
$ \|u^{n+1} -u^n\|_2 \rightarrow 0, \ {\rm as} \ n\rightarrow
\infty. $
\end{enumerate}
\end{lemma}

It can be observed that the properties of $\{u^n\}$ listed in Lemma
{\ref{Lemma_propofnewsequence}} are some direct extensions of those
of $\{x^n\}$. More specifically, Lemma
{\ref{Lemma_propofnewsequence}}(a) can be derived by updating rules
(\ref{updatingindex})-(\ref{updatingrul2}) and the construction
procedure. Lemma {\ref{Lemma_propofnewsequence}}(b) is obtained
directly by the cyclic updating rule. Lemma
{\ref{Lemma_propofnewsequence}}(c) and (d) can be derived by
Property {\ref{Prop_OptimalCondition_n+1}}(b) and the updating rules
(\ref{newupdatingrul1})-(\ref{newupdateFowStep}). Lemma
{\ref{Lemma_propofnewsequence}}(e) can be obtained by Property
{\ref{Prop_SufficientDecrease}} and the definition of $T$
(\ref{def-T}). Lemma {\ref{Lemma_propofnewsequence}}(f) can be
directly derived by Property {\ref{Prop_AsymptoticRegular}}. Besides
Lemma {\ref{Lemma_propofnewsequence}}, the following lemma shows
that the gradient sequence $\{\nabla T(u^n)\}$ satisfies the
so-called relative error condition {\cite{Attouch2013}}, which is
critical to the justification of the convergence of $\{u^k\}$.

\begin{lemma}\label{Lemma_RelativeErrCond}
When $n\geq K-1$, $\nabla T(u^{n+1})$ satisfies
\begin{equation*}
\|\nabla T(u^{n+1})\|_2 \leq b \|u^{n+1} - u^n\|_2,
\end{equation*}
where $b = (\frac{1}{\mu} + K\delta)\sqrt{K},$ with
$$\delta = \max_{i,j= 1,2,\cdots,K} |B_i^TB_{j}|.$$
\end{lemma}

The proof of this lemma is given in Appendix \ref{APPENDIXE}. From
Lemma {\ref{Lemma_propofnewsequence}} (e), the sequence $\{u^n\}$
satisfies the sufficient decrease condition with respect to $T$, and
by Lemma {\ref{Lemma_RelativeErrCond}}, $\{u^n\}$ satisfies the
relative error condition, and also by the continuity of $T$,
$\{u^n\}$ satisfies the so-called continuity condition. Furthermore,
according to {\cite{Attouch2013}} (p. 122), we know that the
function
$$T(u) = \frac{1}{2}\|Bu-y\|_2^2 + \lambda \|u\|_q^q$$
is a {\bf Kurdyka-{\L}ojasiewicz (KL)} function (see Appendix
\ref{APPENDIXG}). Thus, according to Theorem 2.9 in
{\cite{Attouch2013}}, $\{u^n\}$ is convergent.
As a consequence, we can claim the convergence of $\{x^n\}$ as shown
in the following theorem.

\begin{theorem}\label{Thm_MainConv}
Let $\{x^n\}$ be a sequence generated by GAITA with a bounded
initialization. Assume that $0<\mu<L_{\max}^{-1}$, then $\{x^n\}$
converges to a stationary point.
\end{theorem}

According to {\cite{Attouch2013}}, the convergence condition of
JAITA when applied to the ($l_q$LS) problem is $0<\mu <
\|A\|_2^{-2}$. It can be noted that $\max_{i} \|A_i\|_2^2 \leq
\|A\|_2^2$, and hence the condition in Theorem {\ref{Thm_MainConv}}
is generally weaker than that of JAITA. Moreover, as shown by Fig.
{\ref{Weaker_Cond}}, such improvement on the convergence condition
is solid and essential in the sense that there exists a step size
$\mu \in (\|A\|_2^{-2}, L_{\max}^{-1})$ such that JAITA certainly
diverges while GAITA definitely converges with this given step size.

Suppose that $A$ is column-normalized, i.e., $\|A_i\|_2=1$ for any
$i$, then $L_{\max} =1,$ and thus the condition of GAITA becomes
$0<\mu <1.$ In this setting, if further $\mu=1$, then GAITA reduces
to the $l_q$CD algorithm {\cite{lqCD-Marjanovic-2014}} in the
perspective of the iterative form. However, only the subsequential
convergence of the $l_q$CD algorithm can be claimed in
{\cite{lqCD-Marjanovic-2014}} if there is no additional requirement
of $A$. Compared with the $l_q$CD algorithm, there are mainly two
significant improvements. The first one is that we extend the
column-normalized $A$ to a general $A$. Such extension on the model
can improve the flexibility and applicability of GAITA. The second
one, and also the more important one is that the global convergence
of GAITA can be established. It gives a solidly theoretical
guarantee to the use of GAITA.


%
%

\subsection{Convergence to A Local Minimizer}

In this subsection, we mainly answer the second open question
proposed in the end of the subsection 3.1. More specifically, we
will justify that GAITA converges to a local minimizer of the ($l_q$LS) problem under certain conditions.

\begin{theorem}\label{Thm_LocalMin}
Let $\{x^n\}$ be a sequence generated by GAITA with a bounded
initialization. Assume that $0<\mu<L_{\max}^{-1}$, and $x^*$ is the
convergent point of $\{x^n\}$. Let $I = Supp(x^*)$, and
$K=\|x^*\|_0$. Then $x^*$ is a (strictly) local minimizer of
$T_{\lambda}$ if the following condition holds:
\begin{equation}
A_I^TA_I + \lambda q(q-1)\Lambda(x_I^*) \succ 0,
\label{Local_ConvCond}
\end{equation}
where $A_I$ represents the submatrix of $A$ with column restricted
to $I$, $x_I^*$ is the subvector of $x$ restricted to $I$,
$\Lambda(x_I^*) \in \mathbf{R}^{K\times K}$ is a diagonal matrix
with $(|x_i^*|^{q-2})_{i\in I}$ as the diagonal vector, and $M\succ
0$ represents that $M$ is positive definite for any matrix $M$.
\end{theorem}

The proof of this theorem is given in Appendix \ref{APPENDIXF}.
Intuitively, under the condition of Theorem {\ref{Thm_LocalMin}}, it
follows that the principle submatrix of the Henssian matrix of
$T_{\lambda}$ at $x^*$ restricted to the index set $I$ is positive
definite. Moreover, by Lemma {\ref{Lemma_OptimalCondition}} (b), the
following first-order optimality condition holds
\begin{equation*}
A_I^T(Ax^*-y) + \lambda \phi_1(x_I^*) = 0,
\end{equation*}
where $\phi_1(x_I^*) = (qsgn(x^*_{i_1})|x^*_{i_1}|^{q-1},\cdots,
qsgn(x^*_{i_K})|x^*_{i_K}|^{q-1})^T,$ and $i_j \in I, j=1,\cdots,
K.$
These two conditions imply that the second-order optimality
conditions hold at $x^* = (x_I^*,0)$. For any sufficiently small
$h,$ let $x^h = x^* + h,$ then
\begin{align*}
T_{\lambda}(x^h)
&= \frac{1}{2} \|A_I x_I^h -y + A_{I^c}h_{I^c}\|_2^2 + \lambda \sum_{i\in I} |x_i^h|^q +  \lambda \sum_{i\in I^c} |h_i|^q \\
&= \frac{1}{2} \|A_I x_I^h -y \|_2^2 + \lambda \sum_{i\in I} |x_i^h|^q \\
&+ \frac{1}{2} \|A_{I^c}h_{I^c}\|_2^2 + \langle h_{I^c},
A_{I^c}^T(A_I x_I^h -y)\rangle + \lambda \sum_{i\in I^c} |h_i|^q.
\end{align*}
Denote $T_{I^c} = \frac{1}{2} \|A_{I^c}h_{I^c}\|_2^2 + \langle
h_{I^c}, A_{I^c}^T(A_I x_I^h -y)\rangle + \lambda \sum_{i\in I^c}
|h_i|^q.$ Then
\begin{align*}
&T_{\lambda}(x^h) \geq T_{\lambda}(x^*) + T_{I^c}\\
&\geq T_{\lambda}(x^*) + \frac{1}{2} \|A_{I^c}h_{I^c}\|_2^2 +
\sum_{i\in I^c} (\lambda |h_i|^{q} - \|A_{I^c}^T(A_I x_I^h
-y)\|_{\infty} |h_i|),
\end{align*}
where the first inequality holds for the optimality at $x^* = (x_I^*,0)$ and thus, $\frac{1}{2} \|A_I x_I^h -y \|_2^2 + \lambda \sum_{i\in I} |x_i^h|^q \geq T_{\lambda}(x^*)$.
It can be observed that if $h_{I^c}$ is sufficiently small, then the
last part of the above inequality should be nonnegative. Therefore,
$x^*$ should be a local minimizer.

Furthermore, we can drive another two sufficient conditions via
taking advantage of the specific form of the threshold value
(\ref{ThreshValueyLq}). Let $e = \min_{i\in I} |x_i^*|$. Note that
\begin{equation*}
\lambda_{\min}(A_I^TA_I + \lambda q(q-1)\Lambda(x_I^*)) \geq
\lambda_{\min}(A_I^TA_I) + \lambda q(q-1)e^{q-2},
\end{equation*}
where $\lambda_{\min}(M)$ represents the minimal eigenvalue of a
given matrix $M$. Thus, if
\begin{equation}
\lambda_{\min}(A_I^TA_I)>0 \ \text{and}\
0<\lambda<\frac{\lambda_{\min}(A_I^TA_I) e^{2-q}}{q(1-q)},
\label{SuffCond1}
\end{equation}
then the condition of Theorem {\ref{Thm_LocalMin}} holds naturally.

Moreover, by (\ref{ThreshValueyLq}), it holds
\begin{equation}
e\geq \eta_{\mu,q} = (2\lambda\mu (1-q))^{\frac{1}{2-q}}.
\end{equation}
Hence, if $\frac{\lambda_{\min}(A_I^TA_I)}{\max_i \|A_i\|_2^2} >
\frac{q}{2}$ and
\begin{equation}
\frac{q}{2\lambda_{\min}(A_I^TA_I)} < \mu < \frac{1}{\max_{i}
\|A_i\|_2^2}, \label{SuffCond2}
\end{equation}
then the condition (\ref{SuffCond1}) holds and thus
(\ref{Local_ConvCond}) also holds. According to the above analysis,
we can easily obtain the following theorem.

\begin{theorem}\label{Thm_LocalMin1}
Let $\{x^n\}$ be a sequence generated by GAITA with a bounded
initialization. Assume that $0<\mu<L_{\max}^{-1}$, and $x^*$ is the
convergent point of $\{x^n\}$. Let $I = Supp(x^*)$, $K=\|x^*\|_0$,
and $e = \min_{i\in I} |x_i^*|$.
Then $x^*$ is a (strictly) local minimizer of $T_{\lambda}$ if either of the two following conditions satisfies:\\
$ (a)\ \lambda_{\min}(A_I^TA_I)>0,
0<\lambda<\frac{\lambda_{\min}(A_I^TA_I) e^{2-q}}{q(1-q)};
$\\
$ (b)\ \frac{\lambda_{\min}(A_I^TA_I)}{\max_i \|A_i\|_2^2} >
\frac{q}{2}, \frac{q}{2\lambda_{\min}(A_I^TA_I)} < \mu <
\frac{1}{\max_{i} \|A_i\|_2^2}. $
\end{theorem}

Intuitively, the condition (a) in Theorem {\ref{Thm_LocalMin1}}
means that if the smooth part of the ($l_q$LS) problem is strictly
convex and the regularization parameter is sufficiently small, then
the convexity of $T_{\lambda}$ at $x^*$ can be guaranteed by the
convexity of the smooth part. Suppose that $A$ is column-normalized,
i.e., $\|A_i\|_2 =1$ for any $i$, then the condition (b) in Theorem
{\ref{Thm_LocalMin1}} intuitively implies that if the smooth part of
the ($l_q$LS) problem is strongly convex, then the local convexity
of $T_{\lambda}$ at $x^*$ can be guaranteed as long as the step size
$\mu$ is chosen appropriately. Similar conditions are also derived
for the iterative {\em half} thresholding algorithm for solution to
the ($l_q$LS) problem with $q=1/2$ (See Theorems 1 and 2 in
{\cite{ZengHalf2014}}). However, the conditions in this theorem are
a little weaker than those in  {\cite{ZengHalf2014}}.

In {\cite{lqCD-Marjanovic-2014}}, the convergence of the $l_q$CD
algorithm to a local minimizer is justified under a certain scalable
restricted isometry property (SRIP). SRIP is defined as follows.

\begin{definition} \label{SRIP}
{\bf (SRIP {\cite{Beck2011}}).} We say $A$ has the
SRIP($p,\phi,\alpha$) if there exist $\nu_{\phi}, \gamma_{\phi}>0$
satisfying $\gamma_{\phi}/ \nu_{\phi} < \alpha$ such that
$$
\nu_{\phi} \|x\|_2 \leq \|Ax\|_2 \leq \gamma_{\phi} \|x\|_2
$$
holds for every $x \in B_p(\phi):= \{x: \|x\|_p^p \leq \phi\}$, and
$\|\cdot\|_p^p :=\|\cdot\|_0$ for $p=0$.
\end{definition}

Roughly speaking, $\nu_{\phi}$ and $\gamma_{\phi}$ can be viewed as
some type of the minimal and maximal singular values of $A$,
respectively. Thus, SRIP essentially indicates that $A$ possesses a
good condition number. With the definition of SRIP,
{\cite{lqCD-Marjanovic-2014}} demonstrates that if $A$ has the
SRIP($p,\phi,\alpha$) with some $p\geq 0$, then for any $0<q<q^*$
(where $q^* := \min \{1, 2/\alpha^2\}$), the $l_q$CD algorithm
converges to a local minimizer. Particularly, when $\alpha =
\sqrt{2}$, that is, $\gamma_{\phi}/ \nu_{\phi} < \sqrt{2}$, then the
$l_q$CD algorithm converges to a local minimizer for any $0<q<1.$ In
other words, if
\begin{equation}
0<q< \min\{1, \frac{2\nu_{\phi}^2}{\gamma_{\phi}^2}\},
\label{SRIP_Cond}
\end{equation}
then the $l_q$CD algorithm can converge to a local minimizer. It can
be seen from Theorem {\ref{Thm_LocalMin1}} that the condition (b) is
equivalent to
\begin{equation}
0<q<\min \{1, \frac{2\lambda_{\min}(A_I^TA_I)}{\max_i
\|A_i\|_2^2}\}. \label{Local_ConvCond1}
\end{equation}
It is generally hard to compare the conditions (\ref{SRIP_Cond}) and
(\ref{Local_ConvCond1}) directly. However, if $p=0$, then SRIP may
reduce to the standard restricted isometry property (RIP), and in
this case, if further $\phi = K$ (where $K$ is the cardinality of
the support set of $x^*$), then
$$\lambda_{\min}(A_I^TA_I) \geq \nu_K^2\ \text{and}\  \max_i \|A_i\|_2^2 \leq \gamma_K^2.$$
Therefore,
$$\frac{\lambda_{\min}(A_I^TA_I)}{\max_i \|A_i\|_2^2} \geq \frac{\nu_K^2}{\gamma_K^2},$$
which implies that our conditions for convergence to a local
minimizer are generally weaker than that of the $l_q$CD algorithm in
terms of the SRIP.

\section{Numerical Experiments}

In this section, we demonstrate the effects of the algorithmic
parameters on the performance of GAITA. Particularly, we will mainly
focus on the effect of the step size parameter, while the effects of
the regularization parameter $\lambda$ and $q$ can be referred to
{\cite{lqCD-Marjanovic-2014}}. Moreover, a series of experiments are
conducted to show the faster convergence as well as the weaker
convergence condition of GAITA as compared with JAITA.

\subsection{On effect of $\mu$}

For this purpose, we considered the performance of GAITA for the
sparse signal recovery problem, i.e., $y=Ax+\epsilon$, where $x \in
\mathbf{R}^{N}$ was an unknown sparse signal, $A\in \mathbf{R}^{m \times N}$ was the
measurement matrix, $y\in \mathbf{R}^m$ was
the corresponding measurement vector,  $\epsilon$ was the noise and generally $m<N$.
The aim of this problem was to recover the sparse signal $x $ from $y$.
In these experiments, we set $m=250,$ $N=500$ and $k^*
=15,$ where
$k^*$ was the sparsity level of the original sparse signal. The
original sparse signal $x^*$ was generated randomly according to the
standard Gaussian distribution. $A$ was of dimension $m\times N=250
\times 500$ with Gaussian $\mathcal{N}(0,1/250)$ i.i.d. entries and
was preprocessed via column-normalization, i.e., $\|A_i\|_2 =1$ for
any $i$. The observation $y$ was added with 30 dB noise. With these
settings, the convergence condition of GAITA becomes $0<\mu<1.$ To
justify the effect of the step size, we varied $\mu$ from $0$ to $1$,
and considered different $q,$ that is, $q=0.1, 0.3, 0.5, 0.7,
0.9.$ The terminal rule of GAITA was set as the recovery mean
square error (RMSE) $\frac{\|x^{n}-x^*\|_2}{\|x^*\|_2}$  less than a
given precision $tol$ (in this case, $tol = 10^{-2}$).
The regularization parameter
$\lambda$ was set as 0.009 and fixed for all experiments. The
experiment results are shown in Fig. {\ref{Effect_mu}}.

From Fig. {\ref{Effect_mu}}, we can observe that the step size
parameter $\mu$ has almost no influence on the recovery quality of
the proposed algorithm (as shown in Fig. {\ref{Effect_mu}}(a)) while
it significantly affects the time efficiency of the proposed
algorithm (as shown in Fig. {\ref{Effect_mu}}(b)). Basically, we can
claim that the larger step size implies the faster convergence. This
coincides with the common consensus. Therefore, in practice, we
suggest a larger step size like $0.95/L_{\max}$ for GAITA.

\begin{figure}[!t]
\begin{minipage}[b]{.49\linewidth}
\centering
\includegraphics*[scale=0.45]{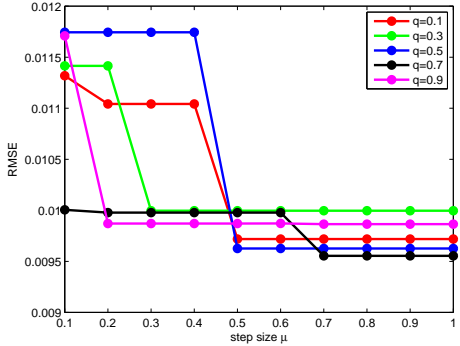}
\centerline{{\small (a) Recovery Error}}
\end{minipage}
\hfill
\begin{minipage}[b]{.49\linewidth}
\centering
\includegraphics*[scale=0.45]{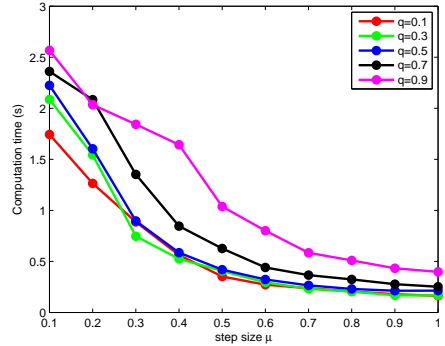}
\centerline{{\small (b) Computational Time}}
\end{minipage}
\hfill \caption{ Experiment for the justification of the effect of
the step size parameter $\mu$ on the performance of GAITA with different
$q$. (a) The trends of recovery error of GAITA with different $q$.
(b) The trends of the computational time of GAITA with different
$q$. } \label{Effect_mu}
\end{figure}

\subsection{Comparison with JAITA}

\subsubsection{Faster Convergence}

We conducted an experiment to demonstrate the faster convergence of
GAITA as compared with JAITA {\cite{L1/2TNN}, \cite{ZengIJT2014}}.
For this purpose, given a sparse signal $x$ with dimension $N=500$
and sparsity $k^*=15,$ shown as in Fig. {\ref{Conv_Rate}}(b), we
considered the signal recovery problem through the observation
$y=Ax,$ where the measurement matrix $A$ and the original sparse
signal $x$ were generated according to the same way in section 4.1.
We then applied GAITA and JAITA  to the ($l_q$LS) problem with two
different $q$, that is, $q=1/2$ and $2/3$, respectively.
In both cases, we took $\lambda=0.001$, $\mu =
\frac{0.95}{\max_{i} \|A_i\|_2^2} (= 0.95)$ for GAITA and
$\mu=0.99\Vert A\Vert_{2}^{-2} (= 0.1676)$ for JAITA. Moreover, the
initial guess was 0. For better comparison, we took every $N$ inner
iterations of GAITA as one iteration since in this case, all
coordinates were updated only once. The experiment results are
reported in Fig. {\ref{Conv_Rate}}.

It can be seen from Fig. {\ref{Conv_Rate}}(a) how the iteration
error ($\Vert x^{n}-x^{\ast}\Vert_{2})$ varies. It can be observed
that GAITA converges much more rapidly than JAITA in both cases. As
shown in Fig. {\ref{Conv_Rate}}(a), the numbers of iterations needed
for GAITA are about 150 in both cases, while much more iterations
are required for JAITA (say, about 1500 and 1700 iterations for
$q=1/2$ and $2/3$, respectively). As justified in
{\cite{ZengHalf2014}, {\cite{ZengIJT2014}}}, JAITA possesses the
eventually linear convergence rate, that is, JAITA will converge
linearly after certain iterations. From Fig. {\ref{Conv_Rate}}(a),
the similar eventually linear convergence rate of GAITA can be
observed. Also, compared with JAITA, much fewer iterations are
required to start such a linear decay. Moreover, Fig.
{\ref{Conv_Rate}}(b) shows that the original sparse signal is
recovered by both GAITA and JAITA with very high accuracies. This
experiment clearly shows the faster convergence as well as
eventually linear convergence rate properties of GAITA.

\begin{figure}[!t]
\begin{minipage}[b]{.49\linewidth}
\centering
\includegraphics*[scale=0.45]{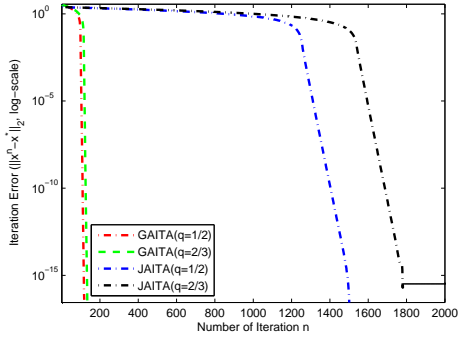}
\centerline{{\small (a) Iteration error}}
\end{minipage}
\hfill
\begin{minipage}[b]{.49\linewidth}
\centering
\includegraphics*[scale=0.45]{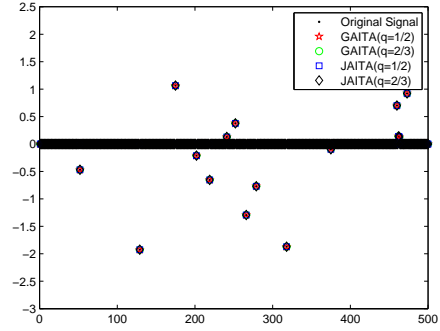}
\centerline{{\small (b) Recovery signal}}
\end{minipage}
\hfill \caption{ Experiment for convergence rate.
(a) The trend of iteration error, i.e., $\|x^{n}-x^*\|_2$. (b)
Recovery signal. The Recovery MSEs of the four cases, that is, GAITA
($q=1/2$), GAITA ($q=2/3$), JAITA ($q=1/2$) and JAITA ($q=2/3$) are
$2.06\times 10^{-8}$, $5.14\times 10^{-9}$, $2.12\times 10^{-8}$ and
$5.28\times 10^{-9}$, respectively. } \label{Conv_Rate}
\end{figure}

\subsubsection{Weaker Condition}

We conducted a set of experiments to demonstrate the convergence
condition of GAITA is weaker than that of JAITA. The experiment
setting was the same as the above subsection. We then applied GAITA
and JAITA to the ($l_q$LS) problem with $q=1/2$. In this setting,
the theoretical condition for convergence of JAITA is $\mu \in
(0,0.1759)$ while the associated condition of GAITA is $\mu \in
(0,1)$. We used different $\mu$ (i.e., $\mu = 0.4, 0.5,0.6,
0.7,0.8,0.9,1$) for both GAITA and JAITA. The figures of the
objective function sequences are shown in Fig. 4.

\begin{figure}[!t]
\begin{minipage}[b]{.49\linewidth}
\centering
\includegraphics*[scale=0.45]{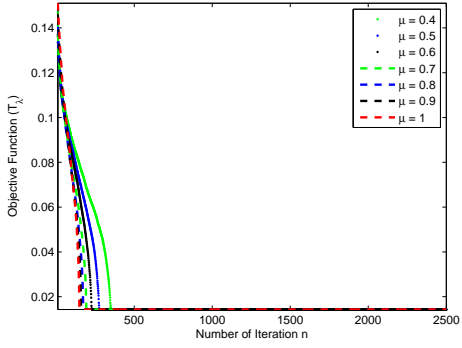}
\centerline{{\small (a) Convergence of GAITA}}
\end{minipage}
\hfill
\begin{minipage}[b]{.49\linewidth}
\centering
\includegraphics*[scale=0.45]{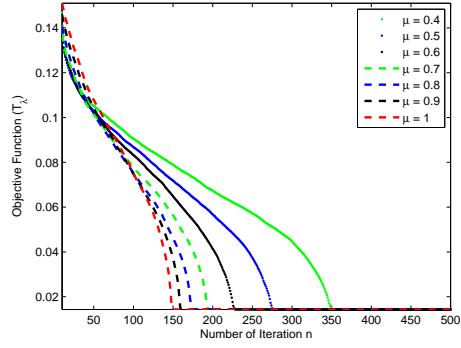}
\centerline{{\small (b) Detail of GAITA}}
\end{minipage}
\hfill
\begin{minipage}[b]{.49\linewidth}
\centering
\includegraphics*[scale=0.45]{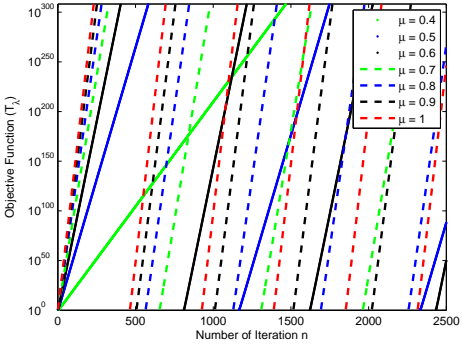}
\centerline{{\small (c) Divergence of JAITA}}
\end{minipage}
\hfill
\begin{minipage}[b]{.49\linewidth}
\centering
\includegraphics*[scale=0.45]{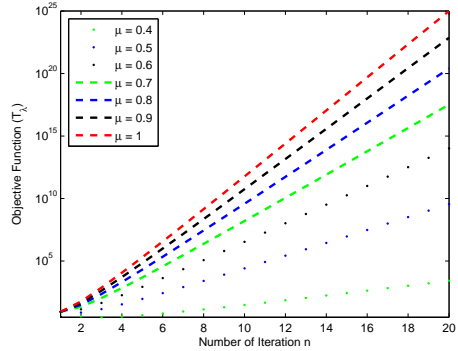}
\centerline{{\small (d) Detail of JAITA}}
\end{minipage}
\hfill \caption{ An experiment that verifies the weaker convergence condition of
GAITA as compared with JAITA.
(a) The trends of the objective function sequences, i.e.,
$\{T_{\lambda}(x^{n})\}$ of GAITA with different $\mu$. (b) The
detail trends of the objective function sequences of GAITA with
different $\mu$. (c) The trends of the objective function sequences
of JAITA with different $\mu$. (d) The detail trends of the objective
function sequences of JAITA with different $\mu$. The regularization
parameter $\lambda$ was taken as 0.001 in all cases. }
\label{Weaker_Cond1}
\end{figure}

From Fig. {\ref{Weaker_Cond1}}, the objective sequences of JAITA
diverge for all $\mu$, while those of GAITA are certainly
convergent. These can be observed detailedly from Fig.
{\ref{Weaker_Cond1}}(b) and (d), respectively. By Fig.
{\ref{Weaker_Cond1}}(a), the objective sequences of GAITA can
converge fast within about 400 iterations in all cases, while those
sequences of JAITA diverge rapidly as shown by Fig.
{\ref{Weaker_Cond1}}(d). When $\mu=1$ and $A$ is column-normalized,
GAITA is reduced to the $l_q$CD method. Fig. {\ref{Weaker_Cond1}}(a)
and (b) show the objective sequence of the $l_q$CD method is
convergent, which can be actually guaranteed by Property
{\ref{Prop_SufficientDecrease}}. It implies that the $l_q$CD method
is weakly convergent as justified in {\cite{lqCD-Marjanovic-2014}}.
However, different from GAITA, the global convergence of the $l_q$CD
method has not been justified if there is no additional condition.

\section{Conclusion}

In this paper, we focused on utilizing the  Gauss-Seidel iteration
rule to the iterative thresholding algorithm for the non-convex
$l_q$ regularized least squares regression problem and developed a
new algorithm called GAITA.
The main contributions of this paper are the
establishment of the convergence of the proposed algorithm. In
summary, we have verified that
\begin{enumerate}
\item[(i)]
GAITA has the finite step convergence of the support set and sign as
long as the step size $0<\mu<1/L_{\max}$. It means that the support
sets and signs of the sequence generated by GAITA can converge and
remain the same within finite iterations.

\item[(ii)]
Under the same condition, the global convergence of GAITA can be
justified. Compared with JAITA like {\em half} algorithm for
$l_{1/2}$ regularization, the convergence condition of GAITA is
weaker than that of JAITA (i.e., $0<\mu<\|A\|_2^{-2}$).

\item[(iii)]
If certain a second-order condition is satisfied at the limit point,
then the limit point can indeed be a local minimizer. Thus, under
these conditions, the proposed algorithm converges to a local
minimizer.

\item[(iv)]
Several numerical experiments are implemented to demonstrate the
effectiveness of GAITA, particularly, the expected faster
convergence and desired weaker convergence condition than JAITA. Also, the
similar eventually linear convergence rate of GAITA can be observed.
However, such rate of convergence property of GAITA has
not been justified in the current paper, and we will study this in
the future work.
\end{enumerate}

%


When it comes to parallel implementation, however, GAITA could have
certain disadvantages because variables that depend on each other
can only be updated sequentially.

\section*{Appendix}

Most of proofs and the description of Kurdyka-{\L}ojasiewicz
inequality are presented in Appendix.

\subsection{Proof of Property
\ref{Prop_SufficientDecrease}}\label{APPENDIXA}

\begin{proof}
Given the current iteration $x^n$, let the coefficient index $i$ be
determined according to (\ref{updatingindex}). According to
(\ref{SVProxOper}) and (\ref{updatingrul1}),
\begin{equation*}
x_i^{n+1} \in \arg \min_{v\in \mathbf{R}}
\left\{\frac{|z_i^n-v|^2}{2} + \lambda\mu |v|^q\right\},
\end{equation*}
where $z_i^n = x_i^n -\mu A_i^T (Ax^n-y)$. Then it implies
\begin{align*}
\frac{1}{2}|\mu A_i^T(Ax^n-y)|^2 + \lambda\mu |x_i^n|^q
\geq \frac{1}{2} |(x_i^{n+1} - x_i^n) + \mu A_i^T (Ax^n -y)|^2 +
\lambda \mu |x_i^{n+1}|^q.
\end{align*}
Some simplifications give
\begin{align}
\lambda |x_i^n|^q - \lambda |x_i^{n+1}|^q
\geq \frac{|x_i^{n+1} - x_i^n|^2}{2\mu} + A_i^T(Ax^n-y)(x_i^{n+1} -
x_i^n). \label{Property2.1}
\end{align}
Moreover, since $x_j^{n+1} = x_j^n$ for any $j\neq i$,
(\ref{Property2.1}) becomes
\begin{align}
\lambda \|x^n\|_q^q - \lambda \|x^{n+1}\|_q^q
\geq \frac{\|x^{n+1} - x^n\|^2}{2\mu} + \langle Ax^n-y, A(x^{n+1} -
x^n) \rangle. \label{Property2.2}
\end{align}
Adding $\frac{1}{2} \|Ax^n -y\|_2^2 - \frac{1}{2} \|Ax^{n+1}
-y\|_2^2$ to both sides of (\ref{Property2.2}) gives
\begin{align}
& T_{\lambda}(x^n) - T_{\lambda}(x^{n+1}) \nonumber\\
& \geq \frac{\|x^{n+1} - x^n\|^2}{2\mu} - \frac{1}{2}\|A(x^n - x^{n+1})\|_2^2 \nonumber\\
& = \frac{\|x^{n+1} - x^n\|^2}{2\mu} - \frac{1}{2}(A_i^TA_i)\|x^n - x^{n+1}\|_2^2 \nonumber\\
& \geq \frac{1}{2}(\frac{1}{\mu} - L_{\max}) \|x^n - x^{n+1}\|_2^2,
\label{Property2.3}
\end{align}
where the first equality holds for
\begin{align*}
\|A(x^n - x^{n+1})\|_2^2
= (A_i^TA_i)|x_i^n - x_i^{n+1}|^2
=(A_i^TA_i)\|x^n - x^{n+1}\|_2^2,
\end{align*}
and the second inequality holds for $A_i^TA_i \leq L_{\max}$.
\end{proof}

\subsection{Proof of Theorem
{\ref{Thm_SubsequentialConv}}}\label{APPENDIXB}

\begin{proof}
By Property {\ref{Prop_SufficientDecrease}}, we know that
$\{T_{\lambda}(x^n)\}$ is a decreasing and lower-bounded sequence,
thus, $\{T_{\lambda}(x^n)\}$ is convergent. Denote the convergent
value of $\{T_{\lambda}(x^n)\}$ as $T^*$. Moreover, by Property
{\ref{Prop_AsymptoticRegular}}, $\{x^n\}$ is bounded, and also by
the continuity of $T_{\lambda}(\cdot)$, there exists a subsequence
of $\{x^n\}$, $\{x^{n_j}\}$ converging to some point $x^*$, which
satisfies $T_{\lambda}(x^*) = T^*$.

Furthermore, by Property {\ref{Prop_AsymptoticRegular}} and
Ostrowski's result (Theorem 26.1, p. 173) {\cite{Ostrowski1973}},
the limit point set $\mathcal{L}$ of the sequence $\{x^n\}$ is
closed and connected.
\end{proof}

\subsection{Proof of Theorem
{\ref{Thm_LimitPoint_StatPoint}}}\label{APPENDIXC}

\begin{proof}
Since the sequence $\{x^n\}$ is bounded, then it has limit points.
Let $x^* \in \mathcal{L}$. We now focus on the $i$-th coefficient of
the sequence with $n=n(i)=jN+i-1$, where $i=1,2,\dots, N$ and
$j=0,1,\dots.$ However, here, we simply use $n$ by which we mean
$n(i)$. Now there exists a subsequence $\{x^{n_1}, x^{n_2},
\cdots\}$ such that
\begin{equation}
\{x^{n_1}, x^{n_2}, \cdots\} \rightarrow x^* \ {\rm and} \
\{x_i^{n_1}, x_i^{n_2}, \cdots\} \rightarrow x_i^*. \label{subseq1}
\end{equation}
Moreover, since the sequence $\{x^{n_1+1}, x^{n_2+1}, \cdots\}$ is
also bounded, thus, it also has limit points. Denoting one of these
by $x^{**}$, then there exists a subsequence $\{x^{l_1+1},
x^{l_2+1}, \cdots\}$ such that
\begin{equation}
\{x^{l_1+1}, x^{l_2+1}, \cdots\} \rightarrow x^{**} \ {\rm and} \
\{x_i^{l_1+1}, x_i^{l_2+1}, \cdots\} \rightarrow x_i^{**},
\label{subseq2}
\end{equation}
where $\{l_1,l_2,\cdots\} \subset \{n_1, n_2, \cdots\}$. In this
case, it holds
\begin{equation}
\{x^{l_1}, x^{l_2}, \cdots\} \rightarrow x^* \ {\rm and} \
\{x_i^{l_1}, x_i^{l_2}, \cdots\} \rightarrow x_i^*, \label{subseq3}
\end{equation}
since it is a subsequence of (\ref{subseq1}). From
(\ref{updateFowStep}) and (\ref{subseq3}), we have
\begin{equation*}
z_i^{l_j} \rightarrow z_i^* \ {\rm as}\ j \rightarrow \infty.
\end{equation*}
Thus, by Lemma {\ref{Lemma_ClosedMap}}, it holds
\begin{equation}
x_i^{**} = \mathcal{T}(z_i^*,x_i^*). \label{Th2.1}
\end{equation}
Moreover, by (\ref{subseq2}), (\ref{subseq3}) and
(\ref{updatingrul2}), it holds
\begin{equation}
x_j^* = x_j^{**} \ {\rm for} \ j\neq i. \label{Th2.2}
\end{equation}

In the following, by the continuity of $T_{\lambda}(\cdot)$ and thus
the continuity of $\Delta(\cdot,\cdot)$ with respect to its
arguments, it holds
\[
\Delta(x^{l_j}, x^{l_j+1}) \rightarrow \Delta(x^*,x^{**}).
\]
Moreover, since the sequence $\{T_{\lambda}(x^n)\}$ is convergent,
then
\[
\Delta(x^{l_j}, x^{l_j+1}) = T_{\lambda}(x^{l_j}) -
T_{\lambda}(x^{l_j+1}) \rightarrow 0 \ {\rm as}\ j\rightarrow
\infty,
\]
which implies
\[
\Delta(x^*,x^{**}) = 0.
\]
Furthermore, by Lemma {\ref{Lemma_DescentFun}}, and
(\ref{Th2.1})-(\ref{Th2.2}), it holds
\begin{equation}
x_i^{**} = x_i^*. \label{Th2.3}
\end{equation}
Combining (\ref{Th2.1}) and (\ref{Th2.3}), we have
\begin{equation}
x_i^* = \mathcal{T}(z_i^*,x_i^*). \label{Th2.4}
\end{equation}

Since $i$ is arbitrary, we have that (\ref{Th2.4}) holds for all $i
\in \{1,\cdots, N\}$. It implies that $x^*$ is a fixed point of
$G_{\mu,\lambda \|\cdot\|_q^q}(\cdot)$, that is, $x^* \in {\mathcal{F}}_q$.
Similarly, since $x^* \in \mathcal{L}$ is also arbitrary, therefore,
$\mathcal{L} \subset {\mathcal{F}}_q$. Consequently, we complete the
proof of this theorem.
\end{proof}

\subsection{Proof of Theorem
{\ref{Thm_SupportConv}}}\label{APPENDIXD}

\begin{proof}
We can note that all the coefficient indices will be updated at
least one time when $n>N$. By Property
{\ref{Prop_OptimalCondition_n+1}}, once the index $i$ is updated at
the $n$-th iteration, then the coefficient $x_i^n$ satisfies:
\[
{\rm either}\  x_i^n=0 \ {\rm or}\  |x_i^n|\geq \eta_{\mu,q}.
\]
Thus, Theorem {\ref{Thm_SupportConv}}(a) holds.

In the following, we prove Theorem {\ref{Thm_SupportConv}}(b) and
(c). By the assumption of Theorem {\ref{Thm_SupportConv}}, there
exits a subsequence $\{x^{n_j}\}$ converges to $x^*$, i.e.,
\begin{equation}
x^{n_j} \rightarrow x^*\ \ \text{as}\ \ j\rightarrow \infty.
\end{equation}
Thus, there exists a sufficiently large positive integer $j_0$ such
that $\|x^{n_j} - x^*\|_2 < \eta_{\mu,q}$ when $j\geq j_0$.
Moreover, by Property {\ref{Prop_AsymptoticRegular}}, there also
exists a sufficiently large positive integer $n^*>N$ such that
$\|x^{n} - x^{n+1}\|_2 < \eta_{\mu,q}$ when $n>n^*$. Without loss of
generality, we let $n^* = n_{j_0}$. In the following, we first prove
that $I^n = I$ and $sgn(x^n) = sgn(x^*)$ whenever $n>n^*$.

In order to prove $I^n=I$, we first show that $I^{n_j} = I$ when
$j\geq j_0$ and then verify that $I^{n+1} = I^n$ when $n>n^*$. We
now prove by contradiction that $I^{n_j} = I$ whenever $j\geq j_0$.
Assume this is not the case, namely, that $I^{n_j} \neq I$. Then we
easily derive a contradiction through distinguishing the following
two possible cases:

\textit{Case 1:} $I^{n_j}\neq I$ and $I^{n_j}\cap I\subset I^{n_j}.$
In this case, then there exists an $i_{n_j}$ such that $i_{n_j}\in
I^{n_j}\setminus I$. By Theorem {\ref{Thm_SupportConv}}(a), it then
implies
\[
\Vert x^{n_j}-x^{\ast}\Vert_{2}\geq|x_{i_{n_j}}^{n_j}|\geq\min_{i\in I^{n_j}%
}|x_{i}^{n_j}|\geq \eta_{\mu,q},
\]
which contradicts to $\Vert x^{n_j}-x^{\ast}\Vert_{2}<\eta_{\mu,q}.$

\textit{Case 2:} $I^{n_j}\neq I$ and $I^{n_j}\cap I=I^{n_j}.$ In
this case, it is obvious that $I^{n_j}\subset I$. Thus, there exists
an $i^{\ast}$ such that $i^{\ast}\in I\setminus I^{n_j}$. By Lemma
{\ref{Lemma_OptimalCondition}}(a), we still have
\[
\Vert
x^{n_j}-x^{\ast}\Vert_{2}\geq|x_{i^{\ast}}^{\ast}|\geq\min_{i\in
I}|x_{i}^{\ast}|\geq \eta_{\mu,q},
\]
and it contradicts to $\Vert
x^{n_j}-x^{\ast}\Vert_{2}<\eta_{\mu,q}$.

Thus we have justified that $I^{n_j}=I$ when $j\geq j_0$. Similarly,
it can be also claimed that $I^{n+1} = I^n$ whenever $n>n^*$.
Therefore, whenever $n>n^*$, it holds $I^n = I$.

As $I^{n}=I$ when $n>n^*$, it suffices to test that
$sgn(x_{i}^{(n)})=sgn(x_{i}^{\ast})$ for any $i\in I$. Similar to
the first part of proof, we will first check that $sgn(x_i^{n_j}) =
sgn(x_i^*)$, and then $sgn(x_i^{n+1}) = sgn(x_i^n)$ for any $i\in I$
by contradiction. We now prove $sgn(x_i^{n_j}) = sgn(x_i^*)$ for any
$i\in I$. Assume this is not the case. Then there exists an
$i^{\ast}\in I$ such that $sgn(x_{i^*}^{n_j})\neq
sgn(x_{i^*}^{\ast})$, and hence,
\[
sgn(x_{i^{\ast}}^{n_j})sgn(x_{i^{\ast}}^{\ast})=-1.
\]
From Lemma {\ref{Lemma_OptimalCondition}}(a) and Theorem
{\ref{Thm_SupportConv}}(a), it then implies
\begin{align*}
\Vert x^{n_j}-x^{\ast}\Vert_{2}
& \geq|x_{i^{\ast}}^{n_j}-x_{i^{\ast}}^{\ast}|=|x_{i^{\ast}}^{n_j}|+|x_{i^{\ast}}^{\ast}|\nonumber\\
& \geq\min_{i\in I}\{|x_{i}^{n_j}|+|x_{i}^{\ast}|\}\geq
2\eta_{\mu,q},
\end{align*}
contradicting again to $\Vert
x^{n_j}-x^{\ast}\Vert_{2}<\eta_{\mu,q}$. This contradiction shows
$sgn(x^{n_j})=sgn(x^{\ast})$. Similarly, we can also show that
$sgn(x^{n+1}) = sgn(x^n)$ whenever $n>n^*$. Therefore, $sgn(x^n) =
sgn(x^*)$ when $n>n^*$.

With this, the proof of Theorem {\ref{Thm_SupportConv}} is
completed.
\end{proof}

\subsection{Proof of Lemma
{\ref{Lemma_RelativeErrCond}}}\label{APPENDIXE}

\begin{proof}
We assume that $n+1 = j^*K+i^*$ for some positive integers $j^* \geq
1$ and $1\leq i^*\leq K$. For simplicity, let
\begin{equation}
i^*=K. \label{i-equalto-K}
\end{equation}
If not, we can renumber the indices of the coordinates such that
(\ref{i-equalto-K}) holds while the iterative sequence $\{u^n\}$
keeps invariant, since the updating rule (\ref{newupdatingindex}) is
cyclic and thus periodic. Such an operation can be described as
follows: for each $n\geq K$, by Lemma
{\ref{Lemma_propofnewsequence}}(b), we know that the coefficients of
$u^n$ are only related to the previous $K-1$ iterates. Thus, we
consider the following a period of the original updating order,
i.e.,
\[
\{i^*+1,\cdots,K,1,\cdots,i^*\},
\]
then we can renumber the above coordinate updating order as
\[
\{1',\cdots,(K-i^*)',(K-i^*+1)',\cdots,K'\},
\]
with
\begin{equation*}
  j'= \left\{
  \begin{array}{cc}
  i^*+j,   & {\rm if}\  1\leq j \leq K-i^*\\
  j-(K-i^*), & {\rm if}\ K-i^*< j \leq K %
  \end{array}%
  \right..
  \label{renum-j}
\end{equation*}%

In the following, we will calculate $\nabla_i T(u^{n+1})$ by a
recursive way for $i=K,K-1,\cdots,1$. Specifically,
\begin{enumerate}
\item[(a)]
For $i=K$, by Lemma {\ref{Lemma_propofnewsequence}}(d), it holds
\begin{equation}
\nabla_K T(u^{n+1}) = (\frac{1}{\mu} - B_K^TB_K) (u_K^n -
u_K^{n+1}). \label{nablaT-n+1-K}
\end{equation}
For any $i=K-1,K-2,\cdots,1,$
\begin{equation*}
\nabla_i T(u^{n+1}) = B_i^T(Bu^{n+1}-y) + \lambda q
sgn(u_i^{n+1})|u_i^{n+1}|^{q-1},
\end{equation*}
and $u_i^{n+1} = u_i^n$. Therefore, for $i=K-1,K-2,\cdots,1,$
\begin{equation}
\nabla_i T(u^{n+1}) = \nabla_i T(u^n) + B_i^T B_K (u_K^{n+1} -
u_K^n). \label{nablaT-red1}
\end{equation}

\item[(b)]
For $i=K-1$, since $n=j^*K+(K-1)$, then by Lemma
{\ref{Lemma_propofnewsequence}}(d) again, it holds
\begin{equation}
\nabla_{K-1} T(u^{n}) = (\frac{1}{\mu} - B_{K-1}^TB_{K-1})
(u_{K-1}^{n-1} - u_{K-1}^{n}).
\end{equation}
By Lemma {\ref{Lemma_propofnewsequence}}(b), it implies
$$
u_{K-1}^{n-1} = u_{K-1}^{n+1}.
$$
Thus,
\begin{equation}
\nabla_{K-1} T(u^{n}) = (\frac{1}{\mu} - B_{K-1}^TB_{K-1})
(u_{K-1}^{n+1} - u_{K-1}^{n}). \label{nablaT-n-K-1}
\end{equation}
Combing (\ref{nablaT-red1}) with (\ref{nablaT-n-K-1}),
\begin{align}
\nabla_{K-1} T(u^{n+1})=
(\frac{1}{\mu} - B_{K-1}^TB_{K-1}) (u_{K-1}^{n+1} - u_{K-1}^{n})
+ B_{K-1}^T B_K (u_K^{n+1} - u_K^n).
 \label{nablaT-n+1-K-1}
\end{align}
Similarly to (\ref{nablaT-red1}), for $i=K-2,K-3,\cdots,1$, we have
\begin{equation}
\nabla_i T(u^{n}) = \nabla_i T(u^{n-1}) + B_i^T B_{K-2} (u_{K-2}^{n}
- u_{K-2}^{n-1}). \label{nablaT-red2}
\end{equation}

\item[(c)]
For any $i=K-j$ with $0\leq j \leq K-1$, by a recursive way, we have
\begin{align}
&\nabla_{K-j} T(u^{n+1}) \nonumber\\
&= \nabla_{K-j} T(u^n) + B_{K-j}^TB_K(u_K^{n+1} -u_K^n)\nonumber\\
&= \nabla_{K-j} T(u^{n-1}) + B_{K-j}^T \sum_{k=0}^{1} B_{K-k}(u_{K-k}^{n+1-k} -u_{K-k}^{n-k})\nonumber\\
&= \cdots \nonumber\\
&= \nabla_{K-j} T(u^{n-j+1})
+ B_{K-j}^T \sum_{k=0}^{j-1} B_{K-k}(u_{K-k}^{n+1-k}
-u_{K-k}^{n-k}). \label{nablaT-n+1-K-j}
\end{align}
Moreover, Lemma {\ref{Lemma_propofnewsequence}}(d) gives
\begin{align}
\nabla_{K-j} T(u^{n-j+1})
= (\frac{1}{\mu} - B_{K-j}^T B_{K-j}) (u_{K-j}^{n-j} -
u_{K-j}^{n-j+1}).
\label{nablaT-n-j+1-K-j}
\end{align}
Plugging (\ref{nablaT-n-j+1-K-j}) into (\ref{nablaT-n+1-K-j}), it
holds
\begin{align}
\nabla_{K-j} T(u^{n+1}) = \frac{1}{\mu} (u_{K-j}^{n-j} - u_{K-j}^{n-j+1})
  + \sum_{k=0}^{j} B_{K-j}^T B_{K-k}(u_{K-k}^{n+1-k}
-u_{K-k}^{n-k}),
\label{nablaT-K-j}
\end{align}
for $j=0,1,\cdots,K-1.$ Furthermore, by Lemma
{\ref{Lemma_propofnewsequence}}(b), it implies
$$
u_{K-k}^{n+1-k} = u_{K-k}^{n+1}
$$
and
$$
u_{K-k}^{n-k} = u_{K-k}^{n}
$$
for $0\leq k \leq K-1$. Therefore, (\ref{nablaT-K-j}) becomes
\begin{align}
\nabla_{K-j} T(u^{n+1}) = \frac{1}{\mu} (u_{K-j}^{n} - u_{K-j}^{n+1})
+ \sum_{k=0}^{j} B_{K-j}^T B_{K-k}(u_{K-k}^{n+1} -u_{K-k}^{n}),
\label{nablaT-K-j*}
\end{align}
for $j=0,1,\cdots,K-1.$
\end{enumerate}

Furthermore, by (\ref{nablaT-K-j*}), it implies
\begin{align}
|\nabla_{K-j} T(u^{n+1})|
&\leq \frac{1}{\mu} |u_{K-j}^{n} - u_{K-j}^{n+1}|
+ \sum_{k=0}^{j} (|B_{K-j}^T B_{K-k}| \cdot |u_{K-k}^{n+1} -u_{K-k}^{n}|)\nonumber\\
& \leq \frac{1}{\mu} |u_{K-j}^{n} - u_{K-j}^{n+1}| + \delta
\|u^{n+1}-u^n\|_1, \label{nablaT-K-j-abs}
\end{align}
for $j=0,1,\cdots,K-1,$ where the second inequality holds for
$$\delta = \max_{i,j=1,\cdots,K} |B_i^T B_j|$$
and
$$\sum_{k=0}^{j} |u_{K-k}^{n+1} -u_{K-k}^{n}| \leq \|u^{n+1}-u^n\|_1.$$
Summing $|\nabla_{K-j} T(u^{n+1})|$ with respect to $j$ gives
\begin{align}
\|\nabla T(u^{n+1})\|_1
&\leq \frac{1}{\mu} \|u^{n+1}-u^n\|_1 + K\delta \|u^{n+1}-u^n\|_1 \nonumber\\
&\leq (\frac{1}{\mu} + K\delta)\sqrt{K}\|u^{n+1}-u^n\|_2,
\label{nablaT-1-norm}
\end{align}
where the second inequality holds for the norm inequality between
1-norm and 2-norm, that is,
\begin{align}
\|u\|_2 \leq \|u\|_1 \leq \sqrt{K} \|u\|_2, \label{NormIneqn}
\end{align}
for any $u\in \mathbf{R}^K$. Also, combining (\ref{NormIneqn}) and
(\ref{nablaT-1-norm}) implies
\begin{equation*}
\|\nabla T(u^{n+1})\|_2 \leq (\frac{1}{\mu} +
K\delta)\sqrt{K}\|u^{n+1}-u^n\|_2.
\end{equation*}
\end{proof}

\subsection{Proof of Theorem {\ref{Thm_LocalMin}}}\label{APPENDIXF}

\begin{proof}
Let $F(x) = \frac{1}{2} \|Ax-y\|_2^2$ and
$$\phi_1(x_I^*) = (qsgn(x^*_{i_1})|x^*_{i_1}|^{q-1},\cdots, qsgn(x^*_{i_K})|x^*_{i_K}|^{q-1})^T,$$
where $i_j \in I, j=1,\cdots, K.$ By Lemma
{\ref{Lemma_OptimalCondition}}(b) we have
\begin{equation}
A_I^T(Ax^*-y) + \lambda \phi_1(x_I^*) = 0. \label{OptCond}
\end{equation}
This together with the condition of the theorem
$$
A_I^TA_I+\lambda q(q-1) \Lambda(x_I^*) \succ 0
$$
imply that the second-order optimality conditions hold at $x^* =
(x_I^*,0)$. For sufficiently small vector $h$, we denote $x_h^* =
(x_I^*+h_I,0)$. It then follows
\begin{equation}
F(x_h^*)+\lambda \sum_{i\in I} |x_i^*+h_i|^q \geq F(x^*) + \lambda
\sum_{i\in I} |x_i^*|^q. \label{SuppPart}
\end{equation}
Furthermore, for any $q\in (0,1)$, it obviously holds that
$$
t^q >(\|\nabla_{I^c} F(x^*)\|_{\infty}+2)t/\lambda,
$$
for sufficiently small $t>0$. By this fact and the differentiability
of $F$, one can observe that for sufficiently small $h$, there hold
\begin{align}
&F(x^*+h)-F(x_h^*) + \lambda \sum_{i\in I^c} |h_i|^q
=\nabla_{I^c}^T F(x^*)h_{I^c} + \lambda \sum_{i\in I^c} |h_i|^q+ o(h_{I^c}) \nonumber\\
&\geq\sum_{i\in I^c} (\|\nabla_{I^c} F(x^*)\|_{\infty} -
[\nabla_{I^c} F(x^*)]_i+1)|h_i| \geq 0. \label{ZeroPart}
\end{align}
Summing up the above two equalities
(\ref{SuppPart})-(\ref{ZeroPart}), one has that for all sufficiently
small $h$,
\begin{equation}
T_{\lambda}(x^*+h) - T_{\lambda}(x^*) \geq 0, \label{LocalMinimizer}
\end{equation}
and hence $x^*$ is a local minimizer.

Actually, we can observe that when $h \neq 0$, then at least one of
the two inequalities (\ref{SuppPart}) and (\ref{ZeroPart}) will hold
strictly, which implies that $x^*$ is a strictly local minimizer.
\end{proof}

\subsection{Kurdyka-{\L}ojasiewicz Inequality}\label{APPENDIXG}

\begin{enumerate}
\item[(a)] The function $f:\mathbf{R} \rightarrow \mathbf{R}\cup \{+\infty\}$ is said to have the Kurdyka-{\L}ojasiewicz property at $x^*\in$ dom $\partial f$ if there exist $\eta \in (0,+\infty]$, a neighborhood $U$ of $x^*$ and a continuous concave function $\varphi:[0,\eta)\rightarrow \mathbf{R}_{+}$ such that:
    \begin{enumerate}
    \item[(i)]
    $\varphi(0) = 0$;

    \item[(ii)]
    $\varphi$ is ${\cal{C}}^1$ on $(0,\eta)$;

    \item[(iii)]
    for all $s\in (0,\eta)$, $\varphi'(s)>0$;

    \item[(iv)]
    for all $x$ in $U\cap \{x: f(x^*)<f(x)<f(x^*) + \eta\}$,
    the Kurdyka-{\L}ojasiewicz inequality holds
    \begin{equation}
    \varphi'(f(x)-f(x^*)) dist(0,\partial f(x)) \geq 1.
    \label{KLIneq}
    \end{equation}

    \end{enumerate}

\item[(b)]
Proper lower semi-continuous functions which satisfy the
Kurdyka-{\L}ojasiewicz inequality at each point of dom $\partial f$
are called KL functions.

\end{enumerate}
Functions satisfying the KL inequality include real analytic
functions, semialgebraic functions and locally strongly convex
functions.

%

\end{document}